\documentclass[12pt,a4paper]{article}
\usepackage[OT1]{fontenc}
\usepackage[latin1]{inputenc}
\usepackage[english]{babel}
\usepackage{amsfonts}
\usepackage{amsthm}
\usepackage{tikz-cd}
\usepackage{amssymb}
\usepackage{mathrsfs}

\usepackage{stmaryrd}
\tikzset{mystyle/.style={sloped, anchor=south}}
\usepackage{mathtools} 

\usepackage{authblk}

\usepackage{float}

\usepackage{amsfonts,amsmath}
\usepackage{fancyhdr,graphicx, xcolor}

\usepackage{subcaption}
\usepackage{t1enc,epsfig,multicol} 
\usepackage{cancel}

\usepackage[symbol]{footmisc}

\def\gr0{ {\text{Gr}(0)} }
\def\gra{ {\text{Gr}(a)} }
\def\Gr{ {\text{Gr}} }
\def\dist{ {\text{dist}} }
\def\length{ {\text{length}} }

\newcommand{\vertiii}[1]{{\left\vert\kern-0.25ex\left\vert\kern-0.25ex\left\vert #1 
		\right\vert\kern-0.25ex\right\vert\kern-0.25ex\right\vert}}

\def\kcal{ {\mathcal K} }
\def\vcal{ {\mathcal V} }
\def\rcal{ {\mathcal R} }
\def\ccal{ {\mathcal C} }
\def\dcal{ {\mathcal D} }

\def\bcal{ {\mathcal B} }

\def\ucal{ {\mathcal{U}} }

\def\hcal{ {\mathcal{H} } }
\def\pcal{ {\mathcal{P}} }

\DeclareMathOperator{\invGr}{invGr}

\def\hgot{ \mathfrak{h} }
\def\hh{ \mathfrak{h} }

\def\pcero{ {\tilde{p}_0} }
\def\px{ {\tilde{p}_{\xb} } }
\def\pp{ {\tilde{p}} }
\def\pz{ {\tilde{p}_{\zb} } }
\def\py{ {\tilde{p}_{\textbf{y}} } }
\def\pcero{ {\tilde{p}_0} }
\def\ptilde{ {\tilde{p}} }
\def\gtilde{ {\tilde{g}} }
\def\qtilde{ {\tilde{q}} }
\def\utilde{ {\tilde{u}} }
\def\vtilde{ {\tilde{v}} }
\def\Xtilde{ {\tilde{X}} }

\def\a{ {\mathcal  A} }
\def\ados{ {{\mathcal {A}}^2} }
\def\lados{{\mathcal{L}}_{\a}(\ados)}
\def\g{\mathcal{G}}

\def\xb{\mathbf{x}}
\def\xbstar{{\mathbf{x}^*}}

\def\zb{\mathbf{z}}

\def\yb{\mathbf{y}}

\def\eb{\mathbf{e}}

\def\B{ {\mathcal  B} }

\def\e{ {\mathcal  E} }

\def\k{ {\mathcal  K} }
\def\f{ {\mathcal  F} }

\def\aa{ {\textbf a} }

\def\CC{ {\mathbb{C}} }
\def\RR{ {\mathbb{R}} }

\DeclareMathOperator{\Asinc}{Asinc}
\DeclareMathOperator{\Exp}{Exp}
\DeclareMathOperator{\Log}{Log}

\DeclareMathOperator{\im}{im}

\DeclareMathOperator{\Grass}{ {\text{Grass}} }
\DeclareMathOperator{\sinc}{sinc}
\DeclareMathOperator*{\diag}{diag}

\DeclareMathOperator*{\ran}{ran}

\newtheorem{teo}{Theorem}[section]
\newtheorem{prop}[teo]{Proposition}
\newtheorem{lem}[teo]{Lemma}
\newtheorem{coro}[teo]{Corollary}
\newtheorem{defi}[teo]{Definition}
\theoremstyle{definition}
\newtheorem{rem}[teo]{Remark}
\newtheorem{ejem}[teo]{Example}
\newtheorem{notacion}[teo]{Notation}

\usepackage[hidelinks]{hyperref} 

\title{The Riemann sphere of a C$^*$-algebra}

\author[1,2]{Esteban Andruchow}
\author[1]{Gustavo Corach}
\author[1]{L\'azaro Recht}
\author[1,2]{Alejandro Varela}
\affil[1]{Instituto Argentino de Matem\'atica ``Alberto P. Calder\'on'', Saavedra 15 3er.�piso, (C1083ACA) Buenos Aires, Argentina}

\affil[2]{Instituto de Ciencias, Universidad Nacional de Gral. Sarmiento, J.	
	M. Gutierrez 1150, (B1613GSX) Los Polvorines, Argentina}
\date{}
\begin{document}

	\maketitle 
	
	\begin{abstract}
		Given the unital C$^*$-algebra $\a$, the unitary orbit of the projector $\pcero=\begin{psmallmatrix}
			1&0\\0&0\end{psmallmatrix}$ in the C$^*$-algebra $M_2(\a)$ of $2\times 2$ matrices with coefficients in $\a$ is called in this paper, \textit{the Riemann sphere $\rcal$ of $\a$}.
			\\
			We show that $\rcal$ is a homogeneous reductive C$^\infty$ manifold of the unitary group $\ucal_2(\a)\subset M_2(\a)$ and carries the differential geometry deduced from this structure (including an invariant Finsler metric). Special attention is paid to the properties of geodesics and the exponential map. If the algebra $\a$ is represented in a Hilbert space $H$, in terms of local charts of $\rcal$, elements of the Riemann sphere may be identified with (graphs of) closed  operators on $H$ (bounded or  unbounded). 
			\\
			In the first part of the paper, we develop several geometric aspects of $\rcal$ including a relation between the exponential map of the reductive connection and the \textit{cross-ratio of subspaces of $H\times H$}. 
			\\
			In the last section we show some applications of the geometry of $\rcal$, to the geometry of operators on a Hilbert space. In particular, we define the notion of \textit{bounded deformation of an unbounded operator} and give some relevant examples.
	\end{abstract}
	


	\tableofcontents 
	\section{Introduction}
		
	This paper is presented (in the spirit of Felix Klein's Erlanger Program) as a sort of ``elliptic'' counterpart of \cite{acr_semipl_poincare_2_capitulos}, where the authors develop aspects of the ``hyperbolic'' Poincar\'e half space of a C$^*$-algebra. Given a unital C$^*$-algebra $\a$ we define the Riemann sphere $\rcal$ of $\a$ as follows. The unitary group $\ucal_2(\a)$ of the C$^*$-algebra $M_2(\a)$ of $2\times 2$ matrices with coefficients in $\a$, operates on the space of projections $\pcal$ of $M_2(\a)$ by inner automorphisms. We call $\rcal$ the $\ucal_2(\a)$-orbit of the projection $\pcero=\begin{psmallmatrix}1&0\\0&0\end{psmallmatrix}$.
	If the algebra $\a$ is faithfully represented in the Hilbert space $H$ and $M_2(\a)$ is correspondingly represented in $H\oplus H$, then $\rcal$ consists of the orthogonal projections in $H\oplus H$ onto subspaces $S$ of the form $S=\utilde\left(H\oplus\{0\}\right)$, $\tilde u\in \ucal_2(\a)$. For example if $\a=B(H)$ and $T:\dcal\to H$ (where $\dcal\subset H$ is dense) is a closed densely defined operator then the orthogonal projection $\ptilde$ onto the graph of $T$ is in $\rcal$ (see Proposition \ref{prop los proy sobre graficos de no acotados estan en R}).
	
	Like in the classical, case where $\a=\mathbb{C}$, we define the ``unitary sphere'' $\kcal\subset \a^2$ as the unitary orbit $\kcal=\{\utilde\eb_1:\utilde\in \ucal_2(\a), \eb_1=\begin{psmallmatrix}		1\\0	\end{psmallmatrix}\}$ 
	and the Hopf fibration $\hgot:\kcal\to\rcal$ by
	$$
	\hgot(\xb)=\xb\xb^*=\begin{pmatrix}
		x_1x_1^*&x_1x_2^*\\x_2x_1^* &x_2x_2^*
	\end{pmatrix}
	$$
where $\xb=\begin{psmallmatrix}
	x_1\\x_2
\end{psmallmatrix}$.
	It is a principal fibration with group $\ucal\subset \a$ as its structure group ($\ucal=$ unitary group of $\a$, acting on $\kcal$ by right multiplication). 

In Section \ref{sec estruct analitica basica} the C$^\infty$ structure of $\rcal$ by means of an appropriate atlas. The principal chart $\mathscr{C}_0=(\vcal_0,\varphi_0,\a)$ of this atlas shows a bijection $\varphi_0$ from $\vcal_0=\{\ptilde\in\rcal: \|\ptilde-\pcero\|<1\}$ onto the algebra $\a$. These projections $\ptilde$ are of the form $\ptilde=\xb\xb^*$ where $\xb=\begin{psmallmatrix}x_1\\x_2\end{psmallmatrix}\in\kcal$ is such that $x_1$ is invertible and the correspondence is given by $\ptilde\mapsto x_2x_1^{-1}$. Like in projective geometry we could call $\begin{psmallmatrix}x_1\\x_2\end{psmallmatrix}$ the ``homogeneous coordinates'' of $\ptilde$ and $a=x_2x_1^{-1}$ the ``affine coordinate''. For this reason we could call this chart the ``projective chart''.

Subsection \ref{The Riemann sphere R as a homogeneous reductive space} is concerned with the geometry of $\rcal$ as a reductive homogeneous space of the group $\ucal_2$. This differential geometric structure contains an affine connection on $\rcal$ and its geodesics. It also contains an invariant Finsler metric which makes $\rcal$ into a metric space where geodesics are minimal curves (see \cite{cpr proy Cstar, cpr proy Banach alg}). The geometry of $\rcal$ as a homogeneous reductive space defines another atlas on $\rcal$ by means of the exponential map of the connection. In Subsection \ref{The inverse of the exponential map in R} we show that the principal chart of this atlas is of the form $(\vcal_0,\Log_\pcero, \mathcal{W})$, where $\vcal_0$ is again $\{\ptilde\in\rcal: \|\ptilde-\pcero\|<1\}$, $\mathcal{W}=\{X\in(T\rcal)_\pcero:\|X\|<\pi/2\}$ and $\Log_{\pcero}$ is the inverse of the exponential map at $\pcero$. We call this chart the ``geodesic chart'' of $\rcal$. 
The relation between the two principal charts has an interesting geometric meaning which we explain in Subsection \ref{secc Geometric interpretation of the logarithm}. Loosely speaking the homogeneous coordinate of $\ptilde$ produce the affine coordinate of $a=x_2x_1^{-1}$ while the geodesic coordinate $X=\begin{psmallmatrix}
	0&a\\a^*&0\end{psmallmatrix}\in (T\rcal)_\pcero$ produces a kind of ``polar coordinate'' of $\ptilde$ given by an ``angle'' and a ``phase'' related to the pair $(\pcero,\ptilde)$.

We devote a final section (Section 5) to the study of examples and applications that we consider relevant. In Subsection 5.3, given an unbounded densely defined closed operator $T$ on a Hilbert space $H$, we show that there exists a unique minimal geodesic on $\rcal$ joining $\pcero$ to $P_{\Gr(T)}$. Notice that $P_{\Gr(T)}$ is in the boundary of $\vcal_0$. In particular we analyze the case of the operator $-i \frac{d}{dx}$ on $L^2[0,1]$. We also study geodesics on $\rcal$ with conjugate points and compute the index of some of these geodesics related to Fredholm operators.
In Section \ref{bounded deformations} we define a notion of (one parameter) bounded deformation of unbounded operators as well as the notion of optimal deformation. The unique minimal geodesic joining $\pcero$ to $P_{\Gr(T)}$, where $T$ is a closed unbounded operator, is an optimal bounded deformation of the unbounded operator $T$.
 In the last section we exhibit types of C$^*$-algebras where $\vcal_0$ is dense in $\rcal$. We remark that $\vcal_0$ is not dense in $\rcal$ when $\a=B(H)$ and $H$ is infinite dimensional.

A second part of this paper will be devoted to the description and uses of a non commutative K\"ahler structure on $\rcal$ which will be defined as an ``elliptic'' counterpart to the one defined in \cite{acr_semipl_poincare_2_capitulos}.

\section{Preliminaries}

We will denote by $\a$ a unital C$^*$-algebra, $\mathcal{G}\subset \a$ its group of invertible elements and $\ucal$ the unitary subgroup of $\mathcal{G}$.
We say that $a\in\a$ is anti-self-adjoint if $a^*=-a$.
The C$^*$-algebra of $2\times 2$ matrices with entries in $\a$ will be denoted by $M_2(\a)$ and the corresponding group of units and unitary subgroup will be denoted by $\mathcal{G}_2$ and $\ucal_2$. Denote by $\a^2$ the right C$^*$ $\a$-module
$$
\a^2=\{\xb=\begin{pmatrix}
	x_1\\x_2
\end{pmatrix}: x_1, x_2\in \a\}
$$
and also write $\a^2_t=\{\hat \xb=\begin{pmatrix}
	x_1 &x_2\end{pmatrix}: x_1,x_2\in \a\}$. We have maps 
\[
\begin{array}{ll}
	\begin{aligned}
	\a^2&\to \a^2_t
	\\
	\xb=\begin{pmatrix}
		x_1\\x_2
	\end{pmatrix}	&\mapsto \xb^*=\begin{pmatrix}
		x_1^* & x_2^*
	\end{pmatrix}	
	\end{aligned}
	&
	\qquad
	\begin{aligned}
		\a^2_t&\to \a^2\\
		\hat\xb=\begin{pmatrix}
			x_1&x_2
		\end{pmatrix}	&\mapsto \hat \xb^*=\begin{pmatrix}
			x_1^* \\ x_2^*
		\end{pmatrix}.
	\end{aligned}
\end{array}
\]
Next we have products
\[
\begin{array}{ll}
	\begin{aligned}
				\a^2&\times\a^2_t \to M_2(\a)
	\\
	\xb\hat \yb&=\begin{pmatrix}
		x_1y_1&x_1y_2\\x_2y_1&x_2y_2
	\end{pmatrix}		
	\end{aligned}
	&
	\qquad
	\begin{aligned}
	 	\a^2_t&\times\a^2 \to  \a
	 	\\
	 	\hat\xb \yb&=	x_1y_1+x_2y_2.
	\end{aligned}
\end{array}
\]

Observe that the inner product in the C$^*$ $\a$-module $\a^2$ is given by $\langle \xb,\yb\rangle=\xb^*\yb=x_1^*y_1+x_2^*y_2$.

	The algebra $M_2(\a)$ is identified with the C$^*$-algebra of $\a$-linear bounded adjointable operators $\mathcal{L}_\a(\ados)$   \cite{lance hilb mod}, where we are fixing the standard basis of $\ados$ given by $\{\eb_1,\eb_2\}$ with $\eb_1=\left(\begin{smallmatrix}
		1\\0
	\end{smallmatrix}
	\right)$, 
	$\eb_2=\left(\begin{smallmatrix}
		0\\1
	\end{smallmatrix}
	\right)$. With this identification, $T\in\lados$,  is represented by the matrix $\tilde{t}=\begin{psmallmatrix}
		t_{11}&t_{12}\\ t_{21}&t_{22}
	\end{psmallmatrix}$:
	$$
	T x=\begin{pmatrix}
		t_{11}&t_{12}\\ t_{21}&t_{22}
	\end{pmatrix}\begin{pmatrix}x_{1}\\ x_{2}
	\end{pmatrix}=
	\begin{pmatrix}
		t_{11}x_1+t_{12}x_2\\ t_{21}x_1+t_{22}x_2
	\end{pmatrix}.
	$$
	Then $T^*$ is given by ${\tilde t}^*=\begin{pmatrix}
		t_{11}^*&t_{21}^*\\ t_{12}^*&t_{22}^*
	\end{pmatrix}.
	$

Note that $\ucal_2$  is the group of $M_2(\a)$ that preserves the quadratic form $(\xb,\yb)\mapsto \xb^*\yb$ when acting on the left by $\xb\mapsto \tilde u \xb$ for $\xb\in\a^2$ and $\tilde u \in \ucal_2$.

\begin{defi}
	A pair of vectors $\xb, \yb\in\a^2$ will be called a \textbf{unitary basis} of $\a^2$ if it is of the form $\xb=\utilde(\eb_1)$  and $\yb=\utilde(\eb_2)$ for $\utilde\in \ucal_2$ where $\eb_1=\begin{psmallmatrix}
		1\\0	\end{psmallmatrix}$ and $\eb_2=\begin{psmallmatrix}
		0\\1	\end{psmallmatrix}$.
\end{defi}
Notice that we have the \textit{Fourier identity} $\zb=\xb \langle\xb,\zb\rangle+\yb\langle\yb,\zb\rangle$ for every $\zb\in\a^2$.

\begin{defi}
	A vector $\xb\in\a^2$ is called a \textit{unitary vector} if it is of the form $\xb=\utilde \eb_1$ for some $\utilde\in\ucal_2$.
\end{defi}
Notice that every unitary vector is the first component of a unitary basis.

Now we come to the central topic of this paper. 
Recall that the unitary group $\ucal_2=\ucal_2(\a)$ operates on the space of all projections of the algebra $M_2(\a)$ by the rule
$$
L_{\utilde}(\ptilde)=\utilde\ptilde\utilde^{-1}\ , \text{ for } \utilde\in\ucal_2.
$$
The geometry related to this action is studied for example in \cite{cpr proy Cstar}.

Let $\pcero$ be the projector $\pcero=\eb_1\eb_1^*=\begin{psmallmatrix}
	1&0\\0&0\end{psmallmatrix}$.

\begin{defi}
	The \textbf{Riemann sphere} $\rcal$ of the C$^*$-algebra $\a$ is the orbit
	$$
	\rcal=\{L_\utilde(\pcero):\utilde\in\ucal_2\}.
	$$
\end{defi}

	A key role in the study of $\rcal$ is played by the space $\kcal$ defined as follows.

\begin{defi}
	%
	We define the \textbf{unitary sphere} $\kcal$ in $\a^2$ as 
	\begin{equation}
		\label{def K}
		\mathcal{K}=\{\xb\in \ados: \exists \tilde{u}\in\ucal_2 \text{ such that }\tilde{u}\eb_1=\xb\}.
	\end{equation}
\end{defi}

\begin{defi}\label{def fibracion Hopf}
	The \textbf{Hopf fibration} over $\rcal$ is the map
	\begin{equation}
		\label{def h fibracion Hopf}
		\mathfrak{h}: \kcal\to \rcal\ ,\ \mathfrak{h}(\xb)=\tilde{p_\xb}
	\end{equation}
	where $\ptilde_\xb=\xb\xb^*$.
	The unitary group $\ucal$ of $\a$ operates by right multiplication on $\kcal$ and is compatible with the projection $\hgot$: $\hgot(\xb u)=\hgot(\xb)$ for $\xb\in\kcal$ and $u\in \ucal$.
\end{defi}

\begin{prop}
	The Hopf fibration is equivariant under the action of $\ucal_2$. More explicitly,
	$$
	\hh(\tilde u\xb)=\tilde u\hh(\xb)\tilde u^*\ ,\ \forall \tilde u\in\ucal_2,\, \xb\in\kcal
	$$
	\begin{proof}
		Indeed, $\hh(\tilde  u\xb)=\tilde u\xb(\tilde u \xb)^*=\tilde u\xb\xbstar\tilde u^*=\tilde u\hh(\xb)\tilde u ^*$.
	\end{proof}
\end{prop}

\begin{notacion}
	Given $\xb=\left(\begin{smallmatrix}
		x_1\\x_2
	\end{smallmatrix}\right)=\tilde{u}\eb_1$ for $\tilde{u}\in \ucal_2$  ($x\in\mathcal{K}$), 
	we will denote with $[\xb]$ the right $\a$-module  generated by $\xb$, that is
	\begin{equation}
		\label{def [x]}
		[\xb]=\xb \a=\{\left(\begin{smallmatrix}
			x_1 a\\x_2 a
		\end{smallmatrix}\right)\in \a^2: a\in \a \}.
	\end{equation}
\end{notacion}

\begin{prop}\label{prop im px es [x]=xA}
	If $\xb\in\mathcal{K}$ and $[\xb]=\xb \a$ is the right $\a$-module  generated by $\xb$, then
	$$
	\im(\tilde{p_\xb})=[\xb].
	$$
\end{prop}
\begin{proof}
	Since $\xb^*\xb=1$ for each $\xb\in \mathcal K$, $\tilde{p_\xb}\xb=(\xb\xb^*)\xb=\xb$ holds and then $\xb\in\im(\tilde{p_\xb})$ which implies that $[\xb]\subset \text{im}(\tilde{p_\xb})$ (since $\text{im}(\tilde{p_\xb})$ is an $\a$-submodule of $\ados$).
	
	The reciprocal is evident. If $\textbf{z}=\tilde{p_\xb}\textbf{w}$ then 
	$\textbf{z}= \xb\xb^*\textbf{w}\in[\xb]$  and hence im$(\tilde{p_\xb})\subset [\xb]$.
\end{proof}

\begin{prop}\label{prop equivs [x] = [z]}
	Suppose that $\xb, \textbf{z}\in \mathcal{K} $. The following statements are equivalent
	\begin{enumerate}
		\item $[\xb]=[\textbf{z}]$
		\item $\tilde{p_\xb}=\tilde{p_\zb}$
		\item $\exists \tilde u\in \ucal$ such that $\textbf{z}=\xb \tilde u$.  
	\end{enumerate}
\end{prop}
\begin{proof}
	(1)$\Leftrightarrow$(2) is evident (see Proposition \ref{prop im px es [x]=xA}). 
	
	(3)$\Rightarrow$(2) is clear after considering that if $u\in\ucal$ and $\textbf{z}=\xb u$, then $\tilde{p_\textbf{z}}=\xb u (\xb u)^*=\xb \xb^*=\tilde{p_\xb}$.
	
	(2)$\Rightarrow$(3)	 If $\tilde{p_\textbf{z}}=\tilde{p_\xb}$, since $\zb=\tilde{p_\zb}\zb=\tilde{p_\xb}\zb=\xb\xb^*\zb$ it is enough to prove that $\xb^*\zb\in\ucal$. That is $(\xb^*\zb)(\zb^*\xb)=1$ and $(\zb^*\xb)(\xb^*\zb)=1$. But $(\xb^*\zb)(\zb^*\xb)=\xb^*\tilde{p_\zb}\xb=\xb^*\tilde{p_\xb}\xb
	=\xb^*\xb=1$ and similarly $(\zb^*\xb)(\xb^*\zb)=1$.
\end{proof}

We mention two natural vector bundles associated to $\rcal$, namely the \textit{tautological vector bundle} $\mathcal{T}\xrightarrow{\text{pr}} \rcal$ (the \textit{bundle of images}) and the \textit{co-tautological vector bundle}  $\mathcal{T}'\xrightarrow{\text{pr'}} \rcal$ (the \textit{bundle of kernels}) defined as follows
\begin{equation}
	\begin{split}\label{eq T y Tprima}
\mathcal{T} &=\{(\ptilde, \xb): \ptilde\in\rcal, \xb\in \im \ptilde\} \text{ and } \text{pr}(\ptilde,\xb)=\ptilde \\
		\mathcal{T'} &=\{(\ptilde, \xb): \ptilde\in\rcal, \xb\in \ker \ptilde\} \text{ and } \text{pr'}(\ptilde,\xb)=\ptilde.
	\end{split}
\end{equation}
Observe that the Hopf fibration is the ``classical bundle of bases'' of the tautological bundle (see Subsection \ref{The canonical connection on the Hopf fibration} for more details). We will also show in  Subsection \ref{The canonical connection on the Hopf fibration} the relation between the co-tautological vector bundle $\mathcal{T}'$ and the tangent bundle $T\rcal$.

\section{The smooth structure of  $\rcal$}\label{sec estruct analitica basica} 

\subsection{The C$^\infty$ structure on  $\kcal$}\label{sec cartas en K} 

Let us start by  constructing a C$^\infty$ structure $\mathfrak{A}_\kcal$ on $\kcal$ (a C$^\infty$ atlas $\mathfrak{A}_\kcal$). 
This is done by identifying open sets in $\kcal$ with appropriate C$^\infty$ manifolds so that the transition functions are C$^\infty$ too.
Define an open neighborhood $\kcal_0$ of $\eb_1$ by
\begin{equation}\label{eq def de K0 dominio carta}
\kcal_0=\{\xb=\left(\begin{smallmatrix}
	x_1\\x_2
\end{smallmatrix}\right)\in\kcal :  x_1\in \g \}
\end{equation}

%

and the map
$$
\psi_0:\kcal_0\to \a\times \ucal \ ,\ \psi_0(\xb)=(x_2 x_1^{-1}, u),
$$
 where $x_1=r u$ is the polar decomposition with $r$ positive and $u$ unitary.

We write now the inverse $\Psi_{0}:\a\times \ucal\to \kcal$ of $\psi_0$
\begin{equation}
	\label{eq Psi0 inversa de psi0}
\Psi_{0}(a,u)=\begin{pmatrix}
	(1+a^*a)^{-1/2} u 
	\\
	a (1+a^*a)^{-1/2} u
\end{pmatrix}.
\end{equation}


We will call the chart given by $\mathcal{C}_0=(\kcal_0,\psi_0,\a\times \ucal)$ the \textit{principal chart} of the C$^\infty$ atlas $\mathfrak{A}_\kcal$.

For each $\tilde u\in \ucal_2$ we describe the chart $\mathcal{C}_\utilde=(\utilde\kcal_0, \psi_\utilde,\a\times \ucal)$ of the atlas $\mathfrak{A}_\kcal$ (by acting with $\utilde$ on the principal chart $\ccal_0$) where $\psi_u=\psi_0\circ \utilde^{-1}:\utilde\kcal_0\to \a\times \ucal$. 
Clearly the atlas $\mathfrak{A}_\kcal=\{\ccal_\utilde:\utilde\in \ucal_2\}$
defines a C$^\infty$ structure on $\kcal$.

\subsection{The C$^\infty$ structure of $\rcal$ 
}\label{secc Differential structure of R}
	Given a unital $C^*$-algebra  $\bcal$ the space
	$$
	\pcal_2=\{p\in \bcal:p^2=p=p^*\}
	$$
	is an C$^\infty$ Banach submanifold of $\bcal$.
	This is well known and details can be found for example in \cite{pr properties of involutions, cpr proy Cstar})
	
	\begin{rem}
	We now recall that the unitary group $\ucal_\bcal$ of $\bcal$ operates on $\pcal$ by inner automorphisms $L_u(p)=upu^*$. This action divides $P$ into orbits and each such orbit is a homogeneous space of the group $\ucal_\bcal$. Moreover, analyzing the infinitesimal situation of this action we can provide each such orbit with a \textit{reductive homogeneous structure}. Details can be found in \cite{cpr proy Cstar}.
	This homogeneous reductive structure provides each orbit with an invariant affine connection and the associated geometry including geodesics, curvature, etc. Details can also be found in \cite{cpr proy Cstar}.
	
	\textit{These ideas apply in our case to the C$^*$-algebra $M_2=M_2(\a)$, and the orbit of $\pcero$ under the action of the unitary group $\ucal_2$, i.e. the Riemann sphere $\rcal$ of the algebra $\a$ and we will use them freely along this paper.
	} 
	\end{rem}
	
	\subsubsection{The C$^\infty$ atlas $\mathfrak{A}_\rcal$ of $\rcal$}\label{secc cartas en K}
	We describe a specific C$^\infty$ atlas on $\rcal$.
	We start by the \textit{principal chart}  $\mathscr{C}_0=(\vcal_0, \varphi_0, \a)$ of this atlas where 
	\begin{equation}
		\label{eq def V0 para p0}
		\vcal_0=\hh\left(\kcal_0\right)=\{\xb \xbstar: \xb\in \kcal , x_1 \text{ invertible }\}\subset \rcal
	\end{equation}
	\text{ and } 
	\begin{equation}
		\label{eq def de fi0 para carta p0}
		\varphi_0:\vcal_0\to \a, \ \varphi_0(\tilde{p})=x_2 x_1^{-1}, \text{ if } \pp =\px \text{ for } \xb\in \kcal_0.
	\end{equation}
%
	
	Observe that if $\pp=\pz$ for another $\zb\in \kcal_0$, then $\zb=\xb u$ for $u\in \ucal$ and $z_2z_1^{-1}=x_2 x_1^{-1}$ and hence $\varphi_0$ is well defined.
	
	Let us verify that $\varphi_0$ is injective. If $\varphi_0(\pp)=\varphi_0(\tilde{q})$ with $\pp=\px$, $\tilde{q}=\py$ for some $\xb, \yb\in \kcal_0$ satisfying $x_2x_1^{-1}=y_2y_1^{-1}$, follows that 
	$ \yb=\left(\begin{smallmatrix}
		y_1\\x_2 x_1^{-1}y_1
	\end{smallmatrix}\right)
	=\left(\begin{smallmatrix}
		1\\x_2 x_1^{-1}
	\end{smallmatrix}\right) y_1$.
	In order to prove that $\pp=\tilde{q}$ it is enough to show that $[\xb]=[\yb]$, since $\im \pp=[\xb]$ and $\im \tilde{q}=[\yb]$. On one hand $\xb=\left(\begin{smallmatrix}
		1\\x_2 x_1^{-1}
	\end{smallmatrix}\right)x_1=\yb y_1^{-1} x_1$ and hence $\xb\in[\yb]$. Analogously, $\yb=\xb x_1^{-1} y_1\in[\xb]$. Therefore $[\xb]=[\yb]$ and $\pp=\tilde q$.
	
	To prove the surjectivity of $\varphi_0$ take any $a\in\a$. We need to find an $\xb\in \kcal_0$ such that $x_2x_1^{-1}=a$, and hence $
	\xb=\left(\begin{smallmatrix}
		x_1\\x_2
	\end{smallmatrix} \right)=\left(\begin{smallmatrix}
		1\\a
	\end{smallmatrix} \right) x_1
	$ 
	should hold.
	To satisfy the condition $\xbstar\xb=1$ we must have that $x_1^* 
	\begin{psmallmatrix}
		1& a^*\end{psmallmatrix}\left(\begin{smallmatrix}
		1\\ a \end{smallmatrix}\right)
	x_1=1
	$. Then $x_1^*(1+a^* a)x_1=1$ which implies that $1+a^*a=(x_1x_1^*)^{-1}$. 
		Every solution of this equation is of the form $x_1=(1+a^*a)^{-1/2} u$, for $u\in\ucal$.
	Now $\xb=\left(\begin{smallmatrix}
		1\\a
	\end{smallmatrix}\right) (1+a^*a)^{-1/2} u$ must satisfy $\xb=\tilde u \eb_1$ for some $\tilde u\in \ucal_2$ and this is the case of
		\begin{equation}
			\label{eq: unitario de U2 para suryectividad}
		\tilde u=	\begin{pmatrix}(1+a^*a)^{-1/2} u & - a^*(1+a a^*)^{-1/2} v
				\\
				a (1+a^*a)^{-1/2} u & (1+a a^*)^{-1/2} v
			\end{pmatrix} \ \text{ , for } v\in\ucal.
		\end{equation}

	We now construct a chart $\mathscr{C}_\utilde=(\vcal_\utilde,\varphi_{\utilde},\a)$ for $\utilde\in\ucal_2$ as follows. We let
	$$
	\vcal_{\utilde}=L_\utilde (\vcal_0) \text{ and } \ 
		\varphi_{\utilde}: \vcal_\utilde\to \a \text{ for }\varphi_\utilde=\varphi_0\circ L_{\utilde^{-1}}.
$$
	Given two charts $\mathscr{C}_\utilde$ and $\mathscr{C}_\vtilde$, $\tilde u, \tilde v\in\ucal_2$ where $\vcal_{\utilde}\cap\vcal_{\vtilde}\neq \emptyset$, let us compute the coordinate change. Let $\xb=\tilde u \eb_1$ and $\yb= \tilde v \eb_1$. We have
	$$
	\left(\varphi_\utilde\circ\varphi_\vtilde^{-1}\right)(a)=\big((\tilde u \tilde v)^*\left(\begin{smallmatrix}
		1\\a
	\end{smallmatrix}\right)\big)_2 \, \Big(\big((\tilde u \tilde v)^*\left(\begin{smallmatrix}
		1\\a
	\end{smallmatrix}\right)\big)_1\Big)^{-1}, \ \forall a\in \varphi_{\vtilde}\left(\vcal_{\utilde}\cap\vcal_{\vtilde}\right).
	$$ 
	If we write $(\tilde u \tilde v)^*=\begin{psmallmatrix}
		c&d\\e&f
	\end{psmallmatrix}\in\ucal_2$  we have
	\begin{equation*}
		\left(\varphi_\px\circ\varphi_\pz^{-1}\right)(a)=(e+f a)(c+d a)^{-1}.
	\end{equation*}
\begin{rem}
Observe that the change of coordinates for two charts in the atlas $\mathfrak{A}_\rcal$ is given by a ``M\"obius'' transformation.
Consequently this atlas defines on $\rcal$ an C$^\infty$ structure.
We shall pursue the study of this complex structure elsewhere.
\end{rem}
	\begin{rem}
		Note that $\tilde p \in \vcal_0$ if and only if there is an element  $\xb=\begin{psmallmatrix}x_1\\x_2 \end{psmallmatrix}\in \im \pcero$ (not necessarily satisfying $\xb^*\xb=1$) such that $x_1$ is invertible, we will call such an $\xb$ a \textit{regular element} in the $\im \pcero$. Also observe that both $\im\pcero$ and $\ker \pcero$ are right $\a$-modules and that every regular element in the $\im \pcero$ is a generator of this $\a$-module.
		Moreover, the correspondence  $\tilde p\in\vcal_0\mapsto a\in \a$ is independent of the choice of the regular element of $\xb\in\im\tilde p$.
		Now choose a faithful representation of $\a$ in a Hilbert space $H$ so that elements $a\in\a$ correspond to operators $a:H\to H$. Consequently $M_2(\a)$ is faithfully represented in $B(H\oplus H)$.
		
	\end{rem}

	\begin{teo} 
		\label{teo sobre proy ortog y esfera}
		 With the previous notations, the following statements are equivalent for $\ptilde\in\rcal$
		\begin{enumerate}
			\item \label{item1teoCarta} $\tilde p\in\vcal_0$ (see \eqref{eq def V0 para p0})
			\item \label{item2teoCarta} $\tilde p$ is the projection $P_{\text{Gr}(a)}\in B(H\oplus H)$ onto the graph of the operator of $a=\varphi(\tilde p)$ (see \eqref{eq def de fi0 para carta p0})
			\item \label{item3teoCarta} $\|\ptilde-\pcero\|<1$
			\item \label{item4teoCarta} $p_{11}$ is invertible if $\ptilde=\begin{pmatrix}
				p_{11}&p_{12}\\ p_{21}&p_{22}
			\end{pmatrix}$.
		\end{enumerate}  
	\end{teo}

	\begin{proof}
		We will denote the elements of $\a$ with the same letters as their representations on $B(H)$.
		First we will prove that item \ref{item1teoCarta} implies \ref{item2teoCarta}. Suppose that $\tilde p=\xb\xb^*=\begin{psmallmatrix}x_1 x_1^*&x_1x_2^*\\x_2x_1^*&x_2 x_2^* \end{psmallmatrix}$ for $\xb\in\kcal$ with $x_1$ invertible. Then $a=\varphi_0(\tilde p)=x_2x_1^{-1}$, $\tilde p\in B(H\oplus H)$ satisfies that $\tilde p^*=\tilde p=\tilde p^2$ and hence it is an orthogonal projection in $B(H\oplus H)$. Moreover, using that $x_1^*x_1+x_2^*x_2=1$ and hence $x_2^*x_2=1-x_1^*x_1$, we can write  $\tilde p\begin{psmallmatrix}		h\\x_2x_1^{-1} h	\end{psmallmatrix}=
		\begin{psmallmatrix}		x_1x_1^*h+x_1x_2^*x_2x_1^{-1}h\\x_2x_1^*h+x_2x_2^*x_2x_1^{-1} h	\end{psmallmatrix}=
		\begin{psmallmatrix}		x_1x_1^*h+x_1(1-x_1^*x_1)x_1^{-1}h\\x_2x_1^*h+x_2(1-x_1^*x_1)x_1^{-1} h	\end{psmallmatrix}=\begin{psmallmatrix}		h\\x_2x_1^{-1} h	\end{psmallmatrix}$. This implies that $Gr(a)\subset \im \tilde p$ for $a=x_2 x_1^{-1}$. Finally, $\tilde p \begin{psmallmatrix}		h\\k	\end{psmallmatrix}=\begin{psmallmatrix}	x_1x_1^*h+x_1x_2^*k\\x_2x_1^*h+x_2x_2^*k	\end{psmallmatrix}
		=\begin{psmallmatrix}	x_1x_1^*h+x_1x_2^*k\\x_2x_1^{-1}(x_1x_1^*h+x_1x_2^*k)	\end{psmallmatrix}
		$ which proves the inclusion $\im \tilde p\subset Gr(a)$.
		
		To prove that item \ref{item2teoCarta} implies \ref{item3teoCarta}, observe first that $\begin{psmallmatrix}h\\k	\end{psmallmatrix}=\begin{psmallmatrix}x\\a x	\end{psmallmatrix} \oplus \begin{psmallmatrix}-a^*y\\y	\end{psmallmatrix}$ for every $\begin{psmallmatrix}h\\k	\end{psmallmatrix}\in H\oplus H$ with
		$\begin{psmallmatrix}x\\a x	\end{psmallmatrix}\in \gra$ orthogonal to $ \begin{psmallmatrix}-a^*y\\y	\end{psmallmatrix}\in\gra^\perp$ (see Lemma \ref{lema Graf T perp}). Then $\|\begin{psmallmatrix}h\\k	\end{psmallmatrix}\|=1$ implies that  $\|\begin{psmallmatrix}x\\a x	\end{psmallmatrix} \|^2 + \|\begin{psmallmatrix}-a^*y\\y	\end{psmallmatrix}\|^2=1$ and hence
		$$
		\|P_\gra-\pcero\|=\sup_{\|\begin{psmallmatrix}h\\k	\end{psmallmatrix}\|=1}
		\|\begin{psmallmatrix}x \\a x	\end{psmallmatrix}- \begin{psmallmatrix}x-a^*y\\0	\end{psmallmatrix}\|
		= \sup_{\|\begin{psmallmatrix}x\\a x	\end{psmallmatrix} \|^2 + \|\begin{psmallmatrix}-a^*y\\y	\end{psmallmatrix}\|^2=1}
		\|\begin{psmallmatrix} a^*y\\ax	\end{psmallmatrix}\|.
		$$
		We know that in general $\|\ptilde-\pcero\|\leq 1$ (see \cite[Corollary 2]{stampfli}). If it es equal to one, there exists a sequence of $\{x_n\}_{n\in\mathbb{N}}$ and $\{y_n\}_{n\in\mathbb{N}}$ such that $\langle a^*y,a^*y\rangle+\langle ax,ax\rangle\to 1$. This would imply that $x_n\to 0$ and $y_n\to 0$ (since $\|x\|^2+\|ax\|^2+\|a^*y\|^2+\|y^2\|=1$), a contradiction. Then 
		$\|P_\gra-\pcero\|<1$.
		
		The proof that item \ref{item3teoCarta} implies \ref{item4teoCarta} follows if we consider the $1,1$ entry of $\ptilde-\pcero$ and observe that $\|\ptilde-\pcero\|<1$. This necessarily implies that $(\ptilde-\pcero)_{1,1}=p_{11}-1$ has norm less than one, and hence $p_{11}$ is invertible.
		
		Finally, if $p_{11}$ is invertible and $\ptilde=\xb\xbstar$ for $\xb\in\kcal$, follows that $p_{11}=x_1 x_1^*$ and therefore $x_1$ is invertible. Then $\ptilde\in\vcal_0$.
	\end{proof}
	
	 \begin{rem}
	 	Theorem \ref{teo sobre proy ortog y esfera} suggests the following consideration. Given an unbounded operator $T:D\to H$ with dense domain $D\subset H$ and closed graph we will show that the orthogonal projection $P_{\text{Gr}(T)}:H\oplus H\to H\oplus H$ belongs to $\rcal(H)$ the Riemann sphere of the algebra $B(H)$. Since we have the obvious embedding $\rcal\subset \rcal(H)$ it is natural to ask about the relative position of $P_{\text{Gr}(T)}$ with respect to $\rcal$. We will give some partial answers to this question in Section \ref{secc The case of unbounded operators}. 
	 	\end{rem}
	 
	 \subsubsection{The tangent map of the principal chart}
	 We describe here the tangent map $(T\varphi_0)_{\tilde p}:(T\rcal)_{\tilde p}\to (T\a)_{\varphi_0(p)} $ ($=\a$) of the principal chart $\varphi_0$ (see \eqref{eq def V0 para p0} and \eqref{eq def de fi0 para carta p0}). Consider the commutative diagram 
	 $$
	 \begin{tikzcd}
	 	\kcal_0\arrow[d, "\hgot|_{\kcal_0}"] \arrow[dr,"\psi_0"]	& \\
	 	\vcal_0 \arrow[r, "\varphi_0"]& \a
	 \end{tikzcd}
	 $$
	 where $\kcal_0=\{\xb\in\kcal: x_1 \text{ is invertible}\}$ (see \ref{secc cartas en K}), $\hgot|_{\kcal_0}$ is the Hopf fibration over $\vcal_0$ and $\psi_0(\xb)=x_2 (x_1)^{-1}$. Now fix $\tilde p\in\rcal$ and $\xb\in\kcal$ such that $\hgot\xb=\tilde p$. Then we have 
	 \begin{equation}\label{ec fla tangente}
	 	(T\varphi_0)_{\tilde p}Y=(T\psi_0)_\xb\, \kappa_\xb(Y), \ \text{ for each } Y\in (T\rcal)_{\tilde p}
	 \end{equation}
	 where $\kappa$ is the structure morphism defined in \ref{subsecc de kapa}.
	 It easy to check that \eqref{ec fla tangente} is independent of the choice of $\xb$ with $\hgot(\xb)=\tilde p$.
	 Explicitly:  $(T\varphi_0)_{\tilde p}Y=(T\psi_0)_\xb Y \xb=(T\psi_0)_\xb \begin{psmallmatrix}
	 	(Y\xb)_1\\(Y\xb)_2
	 \end{psmallmatrix}=(Y\xb)_2 x_1^{-1}- x_2x_1^{-1} (Y\xb)_1 x_1^{-1}$, where we write $Y \xb=\begin{psmallmatrix}(Y\xb)_1\\(Y\xb)_2\end{psmallmatrix}$. 
	 The inverse map $\varphi_0^{-1}$ of $\varphi_0$ is given by 
	 $$
	 \varphi_0^{-1}(a)=\textbf{a} c^2 \textbf{a}^*=\tilde p
	 $$
	 where $\textbf{a}=\begin{psmallmatrix}
	 	1\\a
	 \end{psmallmatrix}\in\a^2$ and $c= (1+a^*a)^{-1/2}\in\a$ (note that $\xb=\textbf{a} c\in \kcal$). The tangent map $(T\varphi_0^{-1})_{a}:(T\a)_a\to (T\rcal)_{\tilde p}$ is given by 
	 $$
	 (T\varphi_0^{-1}) \dot a=\dot \aa c^2\aa^*+\aa c^2 (\dot\aa)^*-\aa b \aa^*
	 $$
	 where $\dot a\in (T\a)_a$ $(=\a)$, $\dot\aa=\begin{psmallmatrix}
	 	0\\ \dot a
	 \end{psmallmatrix}$
	 and $b=c^{-2}(\dot a^*a+a^*\dot a)c^{-2}$.
	 
	 The Finsler structure of $\rcal$ may be translated to a Finsler structure on the manifold $\a$ assigning to the each tangent vector $\dot a\in (T\a)_a$ the norm 
	 $\vertiii{\dot a}=\|\dot \aa c^2\aa^*+\aa c^2 (\dot\aa)^*-\aa b \aa^*\|$ (the standard operator norm of $M_2(\a)$).

	 \subsubsection{Riemann sphere projectors in C$^*$-algebras}\label{subsubsec rsp}
		In this section we give an intrinsic characterization of Riemann spheres of C$^*$-algebras. Let $\mathcal{M}$ be a unital C$^*$-algebra. 
		\begin{defi}
			A self-adjoint projector $p \in\mathcal{M}$ is called a \textbf{Riemann sphere projector (rsp)} if $p$ is conjugated to $1-p$ i.e. there exists an invertible $g\in\mathcal{G}_\mathcal{M}$ such that $g p g^{-1}=1-p$. 
		\end{defi}
		Note that if $p$ is a \textit{rsp} there is a unitary element $u\in\mathcal{M}$ such that $upu^{-1}=1-p$ (this can be shown by an easy argument involving the polar decomposition of $g$ in the definition).
		
		From now on we assume $p$ is a \textit{rsp} in $\mathcal{M}$, $upu^{-1}=1-p$ where $u$ is unitary. Define the subalgebra $\a$ of $\mathcal{M}$ as 
		$$
		\a=p\mathcal{M}p.
		$$
		Note that $\a$ is a C$^*$-algebra with unit $p$. Consider the map $J:\mathcal{M} \to M_2(\a)$ where $J(a)=\begin{psmallmatrix}
			x&y\\z&t
		\end{psmallmatrix}$ where $x=pap$, $y=paup$, $z=pu^{-1}ap$ and $t=pu^{-1}a gu$, with inverse $J^{-1}\begin{psmallmatrix}
		x&y\\z&t
		\end{psmallmatrix}=x+yu^{-1}+uz+utu^{-1}$. 
		
		\begin{prop} $J$ is a C$^*$-algebra isomorphism and $J(p)=\begin{psmallmatrix}
				1&0\\ 0&0
			\end{psmallmatrix}$.
		\end{prop}
		
			The unitary orbit of $p$ in $\mathcal M$ is consequently isomorphic to the Riemann sphere of $M_2(\a)$ and also $p$ and $1-p$ are in the same connected component of the space of projectors of the algebra $\mathcal M$. 
			
			The contents of \ref{subsubsec rsp} are essentially developed in 
		\cite{pr spaces of proj}.

	\subsubsection{The Riemann sphere $\rcal$ as a homogeneous reductive space}\label{The Riemann sphere R as a homogeneous reductive space}
	
	As we have seen the Riemann sphere $\rcal$ of the algebra $\a$ is a subspace $\rcal\subset \pcal_2(\a)$. Also it is clear that the group $\ucal_2$ operates on $\pcal_2(\a)$ by inner automorphisms. In fact $\rcal$ is by definition one orbit of this action, which makes $\rcal$ a homogeneous space of the group $\ucal_2$. This situation is studied in \cite[Section 5]{cpr proy Banach alg} for the general case. In particular $\rcal$ is a homogeneous reductive space of the group $\ucal_2$ and consequently it carries an invariant connection, that we will call the \textit{standard connection}, whose covariant derivative is given by 
		$$
		D_X Y = X\cdot Y + [Y,[X,\tilde p]]
		$$
	where $X\in (T\rcal)_{\tilde p}$, $Y$ is vector field tangent to $\rcal$ in a neighborhood of $\tilde p$ and where $X\cdot Y$ is the directional derivative of $Y$ in the ``ambient'' algebra $M_2(\a)$	
	($X\cdot Y=\frac{d}{dt}\left|_{t=0}\right. Y(\gamma(t))$ for a curve $\gamma(t)$ in $\rcal$, $\gamma(0)=\tilde p$, and $\dot \gamma(0)=X$).
	\begin{rem}
		Given a curve $\ptilde_t\in \rcal$ with $t\in[0,1]$, the differential equation $\frac d{dt} \gtilde_t=[\frac d{dt} \ptilde_t,\ptilde_t]\gtilde_t$ with initial condition $\gtilde_0=1$ has a solution $\gtilde_t\in \ucal_2$ and the action of $\gtilde_t$ on tangent vectors produces the parallel transport of the connection along the curve $\ptilde_t$ (see \cite{cpr proy Banach alg}).
	\end{rem}
	
	Consequently geodesics $\gamma(t)$ in $\rcal$ are defined by the condition 
	$$
	\frac D{dt}\dot\gamma(t)=0
	$$
	and they are explicitly given in the form 
	$$
	\gamma(t)=e^{t\tilde X}\tilde p e^{-t\tilde X}
	$$	
	where $X\in (T\rcal)_{\tilde p}$ and $\tilde X=[X,\tilde p]$. The curve $\gamma(t)$ is the unique geodesic satisfying $\gamma(0)=\tilde p$ and $\dot \gamma(0)=X$.
	Therefore the exponential map is given by $\Exp_\ptilde (X)=\gamma(1)$ where $X\in (T\rcal)_\ptilde$, $\gamma$ is the unique geodesic satisfying $\gamma(0)=\ptilde$ and $\dot \gamma(0)=X$. Observe that $\Exp_{\ptilde}(X)$ is defined for every $X\in (T\rcal)_\ptilde$ and has the explicit form $\Exp_{\ptilde}(X)=e^\Xtilde \ptilde e^{-\Xtilde}$.
\begin{rem}
The exponential map $\Exp_\pcero:(T\rcal)_\pcero\to\rcal$ is bijective from  $V_0=\{X\in (T\rcal)_\pcero: \|X\|<\pi/2\}$ to $U_0=\{\tilde p\in \rcal:\|\ptilde-\pcero\|<1\}$ as we shall see later on (see Theorem \ref{teo Log inversa Exp}).
\end{rem}

%
	
	In what follows we denote by $\sinc$ the analytic function defined by $\sinc(x)=\sin(x)/x$ which is the cardinal $\sin$.
	
	\begin{teo}\label{teo geods desde P(Gr(0)) con CI}
		If $\gamma:[0,1]\to \rcal$ is a geodesic with initial conditions $\gamma(0)=\pcero=\begin{psmallmatrix}
			1&0\\0&0
		\end{psmallmatrix}$ and $\dot \gamma(0)=X=\begin{psmallmatrix}		0&a\\a^*&0		\end{psmallmatrix}$ then, considering the $\im(\pcero)\oplus \im(\pcero)^\perp$ decomposition,
		\begin{equation}
			\label{ec fla geodesicas en teo}
			\begin{split}
				\gamma(t)&=\begin{psmallmatrix} 
					\cos^2|t a^*|& \cos|t a^*|(\sinc|ta^*|)ta\\
					(\sinc|t a|)ta^*\cos|t a^*|& \sin^2|t a^*|\end{psmallmatrix}=
				\begin{psmallmatrix} 
					\cos^2|t a^*|& \sinc\left(2|ta^*|\right)ta\\
					\sinc\left(2|t a|\right)ta^* & \sin^2|t a^*|\end{psmallmatrix}
				\\
				&=
				\begin{psmallmatrix} \cos|t a^*|\\ (\sinc|t a|)ta^*	\end{psmallmatrix}\begin{psmallmatrix} \cos|t a^*| & ta (\sinc|t a|)	\end{psmallmatrix}\\
				&= (\cos^2|t\Xtilde| )\pcero+ (\sin^2|t\Xtilde| )(1-\pcero)+  \sinc (2t|\Xtilde| ) \Xtilde \tilde{\rho_0}
			\end{split}
		\end{equation}
		for $\Xtilde=\begin{psmallmatrix}		0&-a\\a^*&0		\end{psmallmatrix}$ and $\tilde{\rho_0}=\begin{psmallmatrix}		1&0\\0&-1		\end{psmallmatrix}$.
	\end{teo}
	\begin{proof}

		Recall that there is a unique geodesic such that $\dot \gamma(0)=X=\begin{psmallmatrix}		0&a\\a^*&0		\end{psmallmatrix}$ and that can be obtained computing $\gamma(t)=e^{t\tilde X} \pcero e^{-t\tilde X}$ with $\tilde X= \begin{psmallmatrix}		0&-a\\a^*&0		\end{psmallmatrix}$ (see \cite{pr spaces of proj,pr minimality of geod in Grassmann}).

		First, we will describe the unitary given by $e^{\tilde{X}}$ for $\tilde X=\begin{psmallmatrix}		0&-a\\a^*&0		\end{psmallmatrix}$. If we separate the even and odd powers of the series we obtain that
		\begin{equation*}
			\begin{split}
				\tilde X^{2k}&=(-1)^k \begin{psmallmatrix}		|a^*|^{2k}&0\\0&|a|^{2k}		\end{psmallmatrix}, \text{ for } k=0,1,2,\dots
				\ \ \text{ and }\\ 
				\tilde X^{2k+1}&=(-1)^{k+1} \begin{psmallmatrix}	0 &	|a^*|^{2k}a\\-|a|^{2k}a^*&0		\end{psmallmatrix}, \text{ for } k=0,1,2,\dots
			\end{split}
		\end{equation*}
		and then
		\begin{equation*}
			\begin{split}
				\sum_{k=0}^\infty 		\frac{(-1)^k}{(2k)!} \begin{psmallmatrix}		|a^*|^{2k}&0\\0&|a|^{2k}		\end{psmallmatrix}&=\begin{psmallmatrix}		\cos|a^*|&0\\0&\cos|a|		\end{psmallmatrix}
				\ 
				\text{ and }\\
				\
				\sum_{k=0}^\infty 		\frac{(-1)^{k+1}}{(2k+1)!}  \begin{psmallmatrix}	0 &	-|a^*|^{2k}a\\|a|^{2k}a^*&0		\end{psmallmatrix}&=
				\begin{psmallmatrix}		\sinc|a^*|&0\\0&\sinc|a|		\end{psmallmatrix}\begin{psmallmatrix}		0&-a\\a^*&0		\end{psmallmatrix}.
			\end{split}
		\end{equation*}
		This implies that
		\begin{equation}\label{eq e a la Xtilde 1}
			e^{\tilde{X}}=\begin{pmatrix}
				\cos|a^*|& -(\sinc|a^*|)a\\(\sinc|a|)a^*&\cos|a|
			\end{pmatrix}
			=\cos\left|\tilde X\right|+\left(\sinc\left|\tilde X\right|\right)\tilde X
		\end{equation}
		since $|\tilde X|=\begin{psmallmatrix}
			|a^*|&0\\0&|a|
		\end{psmallmatrix}$.
		And for $t\in \mathbb{R}_{\geq 0}$
		\begin{equation}\label{eq e a la t por Xtilde}
			\begin{split}
				e^{t \tilde{X}}&=\begin{pmatrix}
					\cos|t a^*|& -(\sinc|ta^*|)t a\\(\sinc|t a|)t a^*&\cos|t a|
				\end{pmatrix}=\cos|t \tilde X|+\sinc\left(|t\tilde X|\right) \tilde X
			\end{split}.
		\end{equation}

		Then all the geodesics $\gamma$ starting at $\gamma(0)=\pcero$ are of the form
		\begin{equation*}
			\label{ec geodesicas}
			\begin{split}
				\gamma(t)&=e^{t\tilde{X}}\begin{psmallmatrix}1&0\\0&0	\end{psmallmatrix} e^{-t\tilde{X}}
				=\begin{psmallmatrix} 
					\cos^2|t a^*|& \cos|t a^*|(\sinc|ta^*|)ta\\
					(\sinc|t a|)ta^*\cos|t a^*|& \sin^2|t a^*|	\end{psmallmatrix}\\
				&=
				\begin{psmallmatrix} \cos|t a^*|\\ (\sinc|t a|)ta^*	\end{psmallmatrix}\begin{psmallmatrix} \cos|t a^*| & ta(\sinc|t a|)	\end{psmallmatrix}
			\end{split}
		\end{equation*}
		where in the last equality we used that $(\sinc|t a^*|)ta=ta^*\sinc|t a|$.
		
		To obtain the second equality in \eqref{ec fla geodesicas en teo} we can use that $ta^*\cos|t a^*|ta^*\cos|t a^*|=\cos|t a|ta^*$ and that $\cos x\sinc x=\cos x \frac{\sin x}x=\frac12 \frac{\sin(2x)}x=\sinc(2x)$.
		
		The last equality in \eqref{ec fla geodesicas en teo} follows after direct computations.
	\end{proof}
	\begin{rem}
		If the algebra $\a$ is faithfully represented in a Hilbert space $H$ and $\begin{psmallmatrix}
			\xi\\0
		\end{psmallmatrix}\in H\times \{0\}=\im\pcero$, then $e^{t\Xtilde}
		\begin{psmallmatrix}
			\xi\\0
		\end{psmallmatrix}\in \im\gamma(t)$. Observe that $e^{t\Xtilde}
		\begin{psmallmatrix}
		\xi\\0
		\end{psmallmatrix}=
		 \begin{psmallmatrix}
		\cos|t a^*|\xi \\
		 (\sinc|t a|)t a^*\xi
		\end{psmallmatrix}=
		\begin{psmallmatrix}
			\cos|t a^*|\xi \\
			(\sin|t a|)v^*\xi
		\end{psmallmatrix}$ where $a=v|a|$ is the polar decomposition of $a$. For example if $a\geq 0$, we have $e^{t\Xtilde}
		\begin{psmallmatrix}
		\xi\\0
		\end{psmallmatrix}=\begin{psmallmatrix}
		\cos (t a)\xi \\
		\sin(t a)\xi
		\end{psmallmatrix}$, so that in the case where $\xi$ is an eigenvector of $a$, $a\xi=\lambda \xi$, $e^{t\Xtilde}
		\begin{psmallmatrix}
		\xi\\0
		\end{psmallmatrix}$ describes a circular movement in the bidimensional plane generated by $\begin{psmallmatrix}
		\xi\\0
		\end{psmallmatrix}$ and 
		$\begin{psmallmatrix}
		0\\ \xi
		\end{psmallmatrix}$.
	\end{rem}
	\begin{rem}
		Note that if we consider the algebra $\a$ represented in $B(H)$ we can also write the formula \eqref{ec fla geodesicas en teo} as
		\begin{equation*}
			\gamma(t)= 
			\begin{psmallmatrix} 
				\cos^2|t a^*|& \cos|t a^*|(\sin|ta^*|)u\\
				(\sin|t a|)u^*\cos|t a^*|&(\sin|t a|)u^* (\sin|t a^*|) u	\end{psmallmatrix}	
		\end{equation*}
		where $a=u|a|$ is the polar decomposition of $a$ (the partial isometry $u$ might not belong to $\a$).
	\end{rem}

	The space $\rcal$ carries also an invariant Finsler structure given by the
	C$^*$-algebra norm of $M_2(\a)$.
	If $X\in (T\rcal)_{\tilde p}$, $X$ identifies with an element in $M_2(\a)$ and has a corresponding norm.
	This Finsler structure on $\rcal$ allows us to define lengths of curves. In \cite{pr minimality of geod in Grassmann} it is shown that geodesics in $\rcal$ of length less than $\pi/2$ are minimal among curves joining given endpoints.
	
	\subsubsection{The inverse of the exponential map in $\rcal$}\label{The inverse of the exponential map in R}
	
	The standard connection of $\rcal$ defines the exponential map $\Exp_\ptilde:(T\rcal)_\ptilde \to \rcal$ for $\ptilde \in\rcal$. In particular the exponential map $\Exp_\pcero:(T\rcal)_\pcero \to \rcal$ is given by 
	\begin{equation}
		\label{eq formula de Exp}
		\Exp_\pcero(X)=e^{\tilde X}\pcero e^{-\tilde X} , 
	\end{equation}
	where $X\in (T\rcal)_\pcero$ and $\tilde X=[X,\pcero]$ (explicitly, $X=\begin{psmallmatrix}	0&a\\a^*&0	\end{psmallmatrix}$ and $\tilde X=\begin{psmallmatrix}	0&-a\\a^*&0	\end{psmallmatrix}$ for $a\in\a$). It is well known that the exponential map is a diffemorphism of a neighborhood $\mathcal W$ of $0\in(T\rcal)_\pcero$ onto an neighborhood of $\pcero$ in $\rcal$.
	
	In Theorem \ref{teo Log inversa Exp} we will produce an explicit formula for the inverse map $\Log_{\pcero} $ of the exponential map $\Exp_\pcero $. The map $\Log_\pcero$ will be defined on the open set  $U_0=\{\ptilde\in\rcal:\|[\pcero,\ptilde\|<1/2\}$. The mentioned formula involves the real analytic function 
	$$
	\Asinc(x)=\dfrac{\arcsin(x)}{x}, \text{ for } x\in(-1,1).
	$$
	\begin{rem}
		For $\ptilde\in\rcal$, call $\tilde \rho=2\tilde p-1$ (the symmetry associated with $\tilde p$), and observe that the algebra $M_2=M_2^0\oplus M_2^1$ where $M_2^0=\{\tilde a\in M_2: \tilde \rho \tilde a=\tilde a \tilde \rho\}$ and $M_2^1=\{\tilde a\in M_2: \tilde \rho \tilde a=-\tilde a \tilde \rho\}$. Furthermore the mentioned decomposition defines on $M_2$ the structure of $\mathbb{Z}_2$ graded algebra. In this context $(T\rcal)_\ptilde$ may be identified with the self-adjoint part of $M_2^1$. In particular the above formula for $X$ reflects this fact.
		
		At any $\ptilde\in\rcal$ the exponential $\Exp_\ptilde$ is given by
		$
		\Exp_{\ptilde} (X)=e^{\tilde X}\ptilde e^{-\tilde X}$, for  $X$  self-adjoint of degree  $1$  with respect to $\ptilde$  and $\tilde X=[X,\tilde p]$.

	\end{rem}
	 \begin{lem}\label{lem norma corchete p0 p menor 1 medio equiv norma X menor pi sobre 4}
	 	Let $\pcero=\begin{psmallmatrix}	1&0\\0&0\end{psmallmatrix}$ and $\ptilde$ be in $\rcal$ for $\ptilde=e^\Xtilde \pcero e^{-\Xtilde}$ for $\|\Xtilde\|\leq \pi/2$, with $\Xtilde=\begin{psmallmatrix}0&-a\\a^*&0	\end{psmallmatrix}$ and $X=\begin{psmallmatrix}0&-a\\a^*&0	\end{psmallmatrix}\in (T\rcal)_{\pcero}$. Then the following statements are equivalent
	 	\begin{enumerate}
	 		\item $\|[\pcero,\ptilde]\|<1/2$
	 		\item $\|X\|=\|\Xtilde\|<\pi/4$.
	 	\end{enumerate}
	 \end{lem}
	 \begin{proof}
	 First note that $\Xtilde$ anticommutes with $2\pcero-1$ and then 
	 	\begin{equation}\label{eq cuentas norma pcero ptilde}
	 		\begin{split}
	 			1/2&>\|[\pcero,\ptilde]\|=\|\pcero e^{\tilde X} \pcero e^{-\tilde X}-e^{\tilde X} \pcero e^{-\tilde X}\pcero\|
	 			\\
	 			&=\frac 12\|2 \pcero e^{\tilde X} \pcero e^{-\tilde X}-e^{\tilde X} \pcero e^{-\tilde X}+e^{\tilde X} \pcero e^{-\tilde X}-2 e^{\tilde X} \pcero e^{-\tilde X}\pcero)\|
	 			\\
	 			&=\frac 12\|(2 \pcero-1) e^{\tilde X} \pcero e^{-\tilde X}-e^{\tilde X} \pcero e^{-\tilde X}(2 \pcero-1)\|
	 			\\
	 			&=\frac 12\| e^{-\tilde X} \pcero e^{\tilde X}(2 \pcero-1)-e^{\tilde X} \pcero e^{-\tilde X}(2 \pcero-1)\|
	 			\\
	 			&=\frac 12\| e^{-\tilde X} \pcero e^{\tilde X}-e^{\tilde X} \pcero e^{-\tilde X}\|
	 			=\frac 14\| 2e^{-\tilde X} \pcero e^{\tilde X}-1-(2e^{\tilde X} \pcero e^{-\tilde X}-1)\|
	 			\\
	 			&=\frac 14\| e^{-\tilde X} (2\pcero-1) e^{\tilde X}-e^{\tilde X} (2\pcero-1) e^{-\tilde X})\|\\
	 			&=\frac 14\| e^{-2\tilde X}(2\pcero-1) -e^{2\tilde X} (2\pcero-1))\|
	 			=\frac 14 \| e^{-2\tilde X}-e^{2\tilde X} \|=\frac 14 \| 1-e^{4\tilde X} \|.
	 		\end{split}
	 	\end{equation}
	 	Therefore $\|e^{4\Xtilde}-1\|<2$. This implies that $\|e^{4\Xtilde}-1\|<2$ if and only if $\|[\pcero,\ptilde]\|<1/2$. 
	 	But since $\|e^{4\Xtilde}-1\|=\sup_{i\theta \in \sigma(\Xtilde)}|e^{4i\theta}-1|$, 
	 	it follows that $\|[\pcero,\ptilde]\|<1/2$ if and only if $|\theta|<\pi/4$ for $i\theta \in \sigma(\Xtilde)$, which is equivalent to $\|\Xtilde\|<\pi/4$.
	 \end{proof}
	 \begin{lem}\label{lema norma corchete menor que 1 medio entonces p11 mayor o igual que 1 medio}
	 	Let $\pcero$ and $\ptilde\in \rcal$ be such that $\|[\pcero,\ptilde]\|<1/2$, then $p_{11}\geq 1/2$ and $1/2\geq p_{22}$.
	\end{lem}
	\begin{proof}
		The condition \( \|[\pcero,\ptilde]\|<1/2 \) implies that \( \|\Xtilde\|<\pi/2 \), and hence there exists a unique geodesic  
		$
		\gamma(t)=e^{t\Xtilde}\pcero e^{-t \Xtilde}, \quad t\in[0,1],
		$
		between \( \gamma(0)=\pcero \) and \( \gamma(1)=\ptilde \). This condition also implies that \( \|\pcero-\ptilde\|<1 \), and hence \( p_{11} \) is positive definite and invertible.  
		
		Direct calculations give that \( \|[\pcero,\ptilde ]\|=\|p_{12}\|=\|p_{21}\|<1/2 \). Then if \( \ptilde=\xb\xb^* \) for \( \xb=\begin{psmallmatrix}	x_1\\x_2 \end{psmallmatrix}\in\kcal \), 
		where \( x_1=(p_{11})^{1/2} \) is positive definite and invertible, we can write  
		$
		\|x_1x_2^*x_2x_1\|=\|x_1|x_2|^2x_1\| <1/4,
		$  
		since \( p_{21}=x_2x_1 \), \( x_1^2+|x_2|^2=1 \), and \( x_1 \) commutes with \( |x_2| \). Moreover,  we obtain that
		\begin{equation*} 
			\label{eq norma(abs(x2) x1) menor 1/2}
			\| |x_2| x_1\|=\| x_1^{1/2}|x_2| x_1^{1/2}\|<1/2 \text{ which implies } |x_2|x_1\leq 1/2
		\end{equation*}
		since $x_1^{1/2}|x_2| x_1^{1/2}\geq 0$.
		
		We will use now the local cross-section  $\sigma:\{\xb\xb^*:\xb\in\kcal, x_1\in\g\}\to \kcal_0
		$
		of \( \hgot \) (see Theorem \ref{teo fibr hopf estruct analitica} .3 or \eqref{eq seccion local de hgot}). Now define  
		\[
		\begin{psmallmatrix}	\hat x_1(t)\\ \hat x_2(t) \end{psmallmatrix}=\sigma(\gamma(t)).
		\]  
		These entries satisfy that \( \hat x_1(0)=1 \), $\hat x_1(1)=x_1$ and $\hat x_2(1)=x_2$ since \( \gamma(0)=\pcero \) and $\gamma(1)=\ptilde$. Moreover, using the definition of $\sigma$ follows that $x_1(t)>0$ for all $t\in[0,1]$. Now observe that using Lemma \ref{lem norma corchete p0 p menor 1 medio equiv norma X menor pi sobre 4}, since $\|t \Xtilde\|<\pi/4$ for $t\in[0,1]$, 
		then $\|[\pcero,\gamma(t)]\|<1/2$
		and therefore 
		\begin{equation}\label{eq acotacion ||x1(t)|x2(t)| ||<1/2}
			\|x_1(t)|x_2(t)|\|<1/2\ \text{ holds for all } \ t\in[0,1].
		\end{equation}

		The function  $
		g(s)=s\sqrt{1-s^2},  s\in[0,1],
		$
		 is positive in \( (0,1) \), with $g(\sqrt{2}/2)=1/2$ and $g(0)=g(1)=0$.  
		
		Now, suppose that there exists \( t_0\in(0,1] \) such that \( \| \hat x_1(t_0)\|<\sqrt{2}/2 \). Then \( g(\|\hat x_1(t_0)\|)<1/2 \) holds.  
		By the continuity of \( g \), \( \sigma \), and \( \|\cdot\| \), and the fact that \( g(\| \hat x_1(0)\|) = g(1) = 0 \), there exists \( \varepsilon\in(t_0,1) \) such that  
		$
		g(\|\hat x_1(\varepsilon)\|)=1/2.
		$  
		Then,  
		$
		\|x_1(\hat \varepsilon)\|=\sqrt{2}/2,
		$  
		and hence,  
		\[
		\|\hat x_1(\varepsilon)|\hat x_2(\varepsilon)|\|=\|\hat x_1(\varepsilon) (1-\hat x_1(\varepsilon)^2)^{1/2}\|=\|\hat x_1(\varepsilon)\| \sqrt{1-\|\hat x_1(\varepsilon)^2\|}=g(\|\hat x_1(\varepsilon)\|)=1/2.
		\]  
		This contradicts our hypothesis that \( \|[\ptilde,\pcero]\|<1/2 \). The issue arises because we had already established that 
		$\|\hat x_1(t)|\hat x_2(t)|\|<1/2,$
		for all \( t\in [0,1] \) (see \eqref{eq acotacion ||x1(t)|x2(t)| ||<1/2})
		but we reached \( 1/2 \), which is inconsistent with our assumption.
		
		Thus, we conclude that \( \|\hat x_1(t)\|> \sqrt 2/2 \) for all \( t\in [0,1] \), and  
		\[
		\|\hat x_1(t)\|^2=\|\hat x_1(t)^2\|=\|1- |\hat x_2(t)|^2\|=1-\| |\hat x_2(t)|^2\| > 1/2.
		\]  
		This implies that \( \|\hat x_2(t)\|^2< 1/2 \), which gives \( |\hat x_2(t)|\leq  \sqrt{2}/2 \). Then $\sqrt{1-\hat x_1(t)^2}\leq \sqrt{2}/2$ and therefore $1-\hat x_1(t)^2\leq 1/2$ which implies that $1/2\leq \hat x_1(t)^2$ for all $t\in[0,1]$. Hence we have that $1/2\leq \hat x_1(1)^2=x_1^2=\ptilde_{11}$.
		
			In order to prove that $1/2\geq \ptilde_{22}$ observe that since $x_1^2\geq 1/2$ then $1-|x_2|^2\geq 1/2$ and $1/2\geq |x_2|^2$. Thus we obtain that $1/2\geq \||x_2|\|^2 =\|x_2\|=\|x_2^*\|$ and finally that $1/2\geq |x_2^*|^2=p_{22}$.
	\end{proof}
	 
	\begin{teo}\label{teo Log inversa Exp}
			The exponential map $\Exp_\pcero$ (see \eqref{eq formula de Exp}) is a diffeomorphism $\Exp_\pcero:\{X\in (T\rcal)_\pcero: \|X\|<\pi/2\}\to\{\tilde p\in \rcal:\|\ptilde-\pcero\|<1\}$. Moreover, if $U_0=\{\ptilde\in\rcal: \|[\pcero,\ptilde]\|<1/2\}$ and $V_0= \{X\in (T\rcal)_\pcero: \|X\|<\pi/4\}$, there exists an inverse map $\Log_\pcero:U_0\to V_0$ that is a diffeomorphism from $U_0$ to $V_0$ which is given by
		\begin{equation}
			\label{eq formula del Log}
			\Log_\pcero(\ptilde)=\tilde \rho_0\  \Asinc(2|[\pcero,\ptilde]|)\  [\pcero,\ptilde]
		\end{equation}
		where $\tilde{\rho_0}=2\pcero-1$.
		\end{teo}
	
	We call the triple  $(U_0,V_0,\Log_\pcero)$ the \textit{geodesic chart} at $\pcero$. 
	\begin{proof}
		Similar computations to those made in \eqref{eq cuentas norma pcero ptilde} lead to the equivalence between the properties $\|\ptilde-\pcero\|<1$ and $\|\Xtilde\|<\pi/2$ for $\ptilde=e^{\Xtilde}\pcero e^{-\Xtilde}$. 
	Therefore $\Exp_\pcero: \{X\in (T\rcal)_\pcero: \|X\|<\pi/2\}\to\{\tilde p\in \rcal:\|\ptilde-\pcero\|<1\}$ is onto and, since $\|\pcero-\ptilde\|<1$ implies there is a unique geodesic between $\pcero$ and $\ptilde$ of the form $\Exp_\pcero(t\Xtilde)$, then $\Exp_\pcero$ is also injective (see for example \cite[Lemma 2.6]{survey proyectores}). Then $\Exp_\pcero:\{X\in (T\rcal)_\pcero: \|X\|<\pi/2\}\to\{\tilde p\in \rcal:\|\ptilde-\pcero\|<1\}$ is a diffeomorphism (see \cite{pr minimality of geod in Grassmann}).
	
Using again Lemma \ref{lem norma corchete p0 p menor 1 medio equiv norma X menor pi sobre 4} and the fact we mentioned above that $\|\ptilde-\pcero\|<1$ is equivalent to $\|\Xtilde\|<\pi/2$, it can be proved that	$U_0\subset\{\tilde p\in \rcal:\|\ptilde-\pcero\|<1\}$.
	Now we will prove that $\Log_\pcero:U_0\to V_0$ is the inverse of $\Exp_\pcero:V_0\to U_0$.
	
	To prove that the formula \eqref{eq formula del Log} of the inverse holds, we will use the following expression from \eqref{ec fla geodesicas en teo} 
	\begin{equation}\label{fla Exp_pcero ptilde}
		 \Exp_\pcero(X)=\left(\cos^2|\Xtilde|\right)\pcero+\left(\sin^2|\Xtilde|\right)(1-\pcero)+  \sinc\left(2|\Xtilde|\right) \Xtilde \tilde{\rho_0}
	\end{equation}
	for an anti-self-adjoint co-diagonal element $\Xtilde$ such that $X\pcero-\pcero X=\tilde X $.
	Put $X=\tilde{\rho_0} \Asinc(2|[\pcero,\ptilde]|)\,  [\pcero,\ptilde]$. Then $\Exp_\pcero(X)=e^\Xtilde\pcero e^{-\Xtilde}$ for 
	
	\begin{equation}
		\label{eq mod de Xtilde}
		\Xtilde= -\Asinc(2|[\pcero,\ptilde]|)\  [\pcero,\ptilde]\ \text{	 with  }\ 
	|\Xtilde|=\frac12 \arcsin(2|[\pcero,\ptilde]|).
\end{equation}
	Note here that the condition $\|[\pcero,\ptilde]\|=\||[\pcero,\ptilde]|\|<1/2$  implies that   $\arcsin$ and $\Asinc$ are defined and C$^\infty$ in $2|[\pcero,\ptilde]|$. 
	Then, using that $\cos^2\left(\frac{1}{2} \arcsin(2 x)\right)=\frac{1}{2} \left(1+\sqrt{1-4 x^2}\right)$  and that $\ptilde=\xb\xb^*=\begin{psmallmatrix}
	x_{1}x_{1}^*&x_{1}x_{2}^*\\x_{2}x_{1}^*&x_{2}x_{2}^*
\end{psmallmatrix}$ for $\xb\in\vcal_0$ (see Theorem \ref{teo sobre proy ortog y esfera}), we will prove first that $\frac{1}{2} \left(1+\sqrt{1-4 |[\pcero,\ptilde]|^2}\right)_{11}=p_{11}= x_{1}x_{1}^*=x_{1}^2
$ (where we can suppose that $x_{1}$ can be taken invertible and positive). Note that, since $x_{1}^2+|x_{2}|^2=1$ and $x_{1}$ commutes with $|x_{2}|$, we have that 
\begin{equation}
	\begin{split}
\left(2 \cos^2\left(1/2\arcsin\left(2|[\pcero,\ptilde]|\right)_{11}\right)-1\right)^2&=1-4\left(|[\pcero,\ptilde]|^2\right)_{11}=1-4x_{1}^2|x_{2}|^2\\
&=1-4 x_{1}^2(1-x_{1}^2)=(2x_{1}^2-1)^2.
	\end{split}
\end{equation}
We can use here Lemma \ref{lema norma corchete menor que 1 medio entonces p11 mayor o igual que 1 medio} to obtain that $2x_{1}^2-1\geq 0$, since $x_1^2=p_{11}\geq 1/2$.
This implies that $(1-4 x_{1}^2(1-x_{1}^2))^{1/2}=2x_{1}^2-1$ and then $\frac{1}{2} \left(1+\sqrt{1-4 |[\pcero,\ptilde]|^2}\right)_{11}=x_1^2=p_{11}$ which is the equality $\Exp_\pcero(X)_{11}=\ptilde_{11}$.

For the $\ptilde_{22}$ entry we can reason similarly, but using that $\sin^2\left(\frac{1}{2} \arcsin(2 x)\right)=\frac{1}{2} \left(1-\sqrt{1-4 x^2}\right)$, to obtain that
\begin{equation*}
	\begin{split}
		(1-2\sin^2(	1/2 \arcsin(2|[\pcero,\ptilde]|_{22})))^2		&=
		1-4(|[\pcero,\ptilde]|^2)_{22}
		=1-4|p_{12}|^2\\
		&=1-4(x_1x_2^*)^*x_1x_2^*=1-4x_2 x_1^2x_2^*\\
		&=1-4x_2(1-x_2^*x_2)x_2^*=1-4(|x_2^*|^2-|x_2^*|^4)\\
		&=(1-2|x_2^*|^2)^2.
	\end{split}
\end{equation*}
Then, since $|x_2^*|^2=p_{22}\leq 1/2$ (see Lemma \ref{lema norma corchete menor que 1 medio entonces p11 mayor o igual que 1 medio}) $1-2|x_2^*|^2\geq 0$ holds, and hence we obtain that $1-2\sin^2(	1/2 \arcsin(2|[\pcero,\ptilde]|_{22})=1-2|x_2^*|^2$. Therefore 
$$
\Exp_\pcero(X)_{22}=\sin^2(	1/2 \arcsin(2|[\pcero,\ptilde]|_{22})=|x_2^*|^2=p_{22}.
$$

Considering the last term of \eqref{fla Exp_pcero ptilde}, the codiagonal of $\Exp_\pcero(X)$, observe that if $[\pcero,\ptilde]=|[\pcero,\ptilde]|\tilde\nu$ is the polar decomposition of $[\pcero,\ptilde]$, we can write (see \eqref{eq mod de Xtilde})
\begin{equation*}
	\begin{split}
\sinc(2|\Xtilde|)\Xtilde\tilde\rho_0&=-\frac12\sin(2|\Xtilde|)\tilde\nu \tilde\rho
=-\frac12\sin(\arcsin(2|[\pcero,\ptilde]|))\tilde\nu\tilde\rho_0=-|[\pcero,\ptilde]|\tilde\nu\tilde\rho_0\\
&=-[\pcero,\ptilde]\tilde\rho_0.
	\end{split}
\end{equation*}

Then, since  $[\pcero,\ptilde]=\begin{psmallmatrix}	0&p_{12}\\-p_{21}&0\end{psmallmatrix}$ 
the codiagonal of $\ptilde$ coincides with $-[\pcero,\ptilde]\rho_0$.

Therefore we have proved that $\Exp_\pcero(-\tilde{\rho_0} \Asinc(2|[\pcero,\ptilde]|)\  [\pcero,\ptilde])=\ptilde$, for $\ptilde\in U_0$.
	\end{proof} 

	\begin{rem}
		The formula \eqref{eq formula del Log} does not hold as the inverse of $\Exp_\pcero$ in the domain $\{\ptilde\in\rcal: \|\ptilde-\pcero\|<1\}\supset \{\ptilde\in\rcal: \|[\pcero,
		\ptilde\|<1/2\}$. An example where $\|\ptilde-\pcero\|<1$ but $\|[\pcero,\ptilde]\|>1/2$ hold and $\Exp_\pcero\left(\tilde \rho_0\  \Asinc(2|[\pcero,\ptilde]|)\  [\pcero,\ptilde]\right)(\ptilde)\neq \ptilde$, is 
		$\ptilde =	\begin{psmallmatrix}
					\cos^2(\pi/3) & \sin (\pi/3 ) \cos (\pi/3 ) \\		
					\sin (\pi/3) \cos (\pi/3) & \sin^2(\pi/3) 
		\end{psmallmatrix}=
			\begin{psmallmatrix}
			 {1}/{4} & {\sqrt{3}}/{4} \\
			{\sqrt{3}}/{4} & {3}/{4} 
			\end{psmallmatrix}\in M_2(\CC).$
\end{rem}
	
	\begin{rem} 
		In general, if $\Exp_\pcero(X)=\ptilde$ and $\|X\|< \pi/2$, we will say that $X\in (T\rcal)_\pcero$ is the \textit{geodesic coordinate} of $\ptilde$. In this way we have geodesic coordinates in $(T\rcal)_\pcero$ for points $\ptilde\in\rcal$ such that $\|[ \pcero,\ptilde]\|<1/2$.
	\end{rem}
	
	\begin{rem}
	Consider a representation of the algebra $\a$ into a Hilbert space $H$ and the corresponding representation of $M_2(\a)$ in $H\oplus H$. Next write $[\pcero,\ptilde]=|[\pcero,\ptilde]| \utilde$,  the polar decomposition of $[\pcero,\ptilde]$, where $\utilde$ is the partial isometry and observe that $\utilde$ commutes with $|[\pcero,\ptilde]|$ since $[\pcero,\ptilde]$ is anti self-adjoint. 
	Then we can write 
	$$
	\Log_\pcero(\ptilde)=  \frac12 \arcsin(2|[\pcero,\ptilde]|)  \tilde \rho_0 \utilde
	$$
	the polar decomposition of $\Log_\pcero(\ptilde)$.
	In this formula we may interpret the positive part $\frac12 \arcsin(2|[\pcero,\ptilde]|)$ as a kind of ``unoriented'' angle between $\pcero$ and $\ptilde$ and the partial isometry $\tilde{\rho_0} \utilde$ as a \textit{partial imaginary unit} in the sense that $(\tilde{\rho_0}\utilde)^2=-\tilde q$ where $\tilde q$ is the projection $\utilde^*\utilde$.
	
\end{rem}

\subsubsection{Geometric interpretation of the logarithm}\label{secc Geometric interpretation of the logarithm}
	We start with an example. 
	
	Let $\a=\CC$ so $M_2=M(2,\CC)$ is the C$^*$-algebra of $2\times 2$ complex matrices. The $\ucal_2$ orbit of $\textbf{e}_1=\begin{psmallmatrix}1\\0\end{psmallmatrix}$  is $\kcal=S^3$ the unit sphere of $\CC^2$. Also note that the Riemann sphere $\rcal$ is the original Riemann sphere which is here represented by the orbit of the projector $\pcero=\begin{psmallmatrix}
		1&0\\0&0\end{psmallmatrix}$ under the action of $\ucal_2$ (of course $\rcal$ is diffeomorphic to the projective line $\mathbb{P}^1(\CC)``="S^2$). Finally the Hopf fibration is the original Hopf fibration given by $\hgot:\kcal\to\rcal$
		$$\zb\mapsto \hgot(\zb)=\ptilde_z, \ \text{ where }
		\zb=\begin{pmatrix}
			z_1\\z_2\end{pmatrix} \ ,\ \tilde p_z=\zb\zb^*= \begin{pmatrix}
			|z_1|^2&z_1\overline{z_2}\\ \overline{z_1}z_2&|z_2|^2
			\end{pmatrix}.
	$$

	Let $\ptilde\in\rcal$ be such that $\|[\pcero,\ptilde]\|<1/2$ and write $X=\Log_\pcero(\ptilde)$. 
	According to the explicit formula for $\Log_\pcero$ we may write 
	$$
	X= \frac{\arcsin(2|[\pcero,\ptilde]|)}2  \tilde{\rho_0} \utilde 
	$$
	where $[\pcero,\ptilde]=|[\pcero,\ptilde]| \utilde$ ($\utilde$ partial isometry) is the polar decomposition  of  $[\pcero,\ptilde]$.
	Notice that $|[\pcero,\ptilde]|$ commutes with $[\pcero,\ptilde]$.
	Observe that $|[\pcero,\ptilde]|= |z_1|\, |z_2|$ (a scalar in $M(2,\CC)$).
	Now since $ |z_1|^2+ |z_2|^2=1$ there is a unique angle $0\leq \varphi\leq \pi/2$ such that $|z_2|=\sin(\varphi)$ and  $|z_1|=\cos(\varphi)$. Therefore the positive part of $X$ is exactly $\varphi$, so 
	$$
	X=\varphi\ \tilde{\rho_0} \utilde
	$$ 
	in its polar decomposition.
	Finally, the positive part of the logarithm of $\ptilde$ is the Finsler distance $\text{dist}(\pcero,\ptilde)$ in the Riemann sphere.
	
\begin{figure}[h]
	\centering{\includegraphics[width=0.5\textwidth]{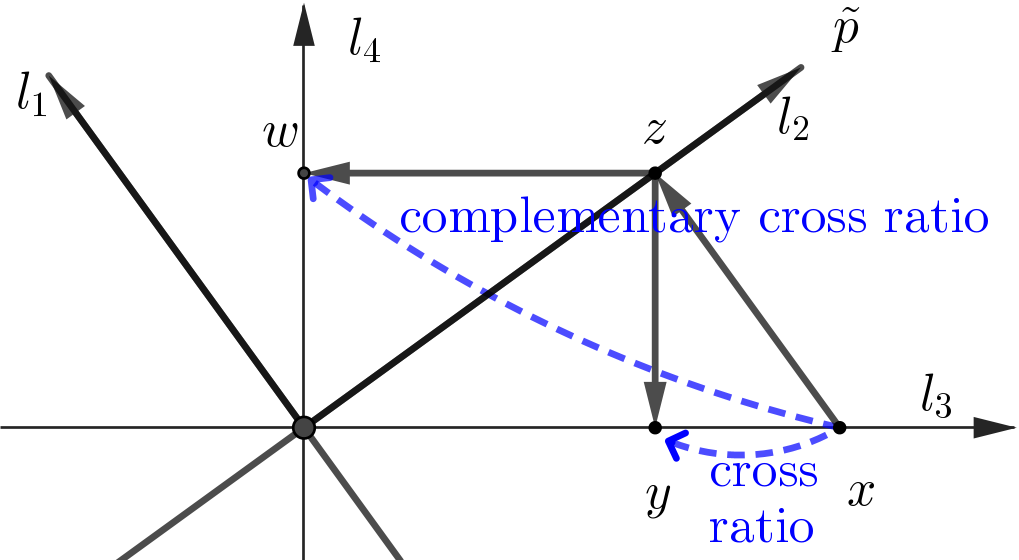} }
		\caption[short]{Cross ratio and complementary cross ratio.}
		\label{razon doble}
		\end{figure}

	Next we produce a geometric interpretation of $[\pcero,\ptilde]$ ($=|[\pcero,\ptilde]|\utilde$). Consider the Figure \ref{razon doble}.
	In it we see schematically $\pcero$, $\ptilde$ and the (complex) lines $l_1=\ker(\ptilde)$, $l_2=\im(\ptilde)$, $l_3=\im(\pcero)$ and $l_4=\ker(\pcero)$. The correspondence $x\in l_3$ maps to $y\in l_3$ defines a linear map $y=\alpha x$ from $l_3$ to $l_3$. The number $\alpha$ is the classical cross ratio of the ordered four points $l_1$, $l_2$, $l_3$, $l_4$ in the complex projective line. In our case $[\pcero,\ptilde]$ has the form 
	$$
	[\pcero,\ptilde]=
	\begin{pmatrix}
		0&z_1\overline{z_2}\\ -\overline{z_1}z_2&0
	\end{pmatrix}
	$$
	and therefore $[\pcero,\ptilde]$ maps $l_3$ into $l_4$ and $l_4$ into $l_3$.
	We only describe the map $l_3\to l_4$ (the other one is similar).
	The correspondence $x\to w$ in the picture has a matrix $\beta=-\overline{z_1}z_2$ which determines $[\pcero,\ptilde]$. We call this geometric construction the complementary cross ratio.
	
	\bigskip
	
	With this example in mind we now turn to the general case.
Observe first that the inverse (see \cite{pr minimality of geod in Grassmann}) of the diffeomorphism $\Exp_\pcero:\{\ptilde\in\rcal: \{X\in(T\rcal)_\pcero:\|X\|<\pi/2\}\to \|\ptilde-\pcero\|<1\}$ defined in \eqref{eq formula de Exp}, allows us to determine an angle between $\ptilde$ and $\pcero$ using the polar decomposition of the corresponding $X\in(T\rcal)_\pcero$ in some representation of $\a$.
	In what follows, we consider the formula of $\Log_{\pcero}$ from \eqref{eq formula del Log} to obtain an expression of this angle.
	Let $\ptilde\in \rcal$ be such that $\|[\pcero,\ptilde]\|<1/2$. Then we have that 
	$ \Log_{\pcero} \ptilde=\Asinc(2|[\pcero,\ptilde]|) (\tilde{\rho_0} [\pcero,\ptilde])$ (see Theorem \ref{teo Log inversa Exp}).
	
We can represent the algebra $\a$, and correspondingly $M_2(\a)$, faithfully in a Hilbert space  $H$ (resp. $H\times H$) and refer the polar decompositions to this representation.  
Write the right polar decomposition of the bracket $[\pcero,\ptilde]$ as $[\pcero,\ptilde]=|[\pcero,\ptilde]| \vtilde$, where $\vtilde$ is the partial isometry. Note that $[\pcero,\ptilde]=(-\vtilde^*)|[\pcero,\ptilde]|$ is the left polar decomposition.

	We claim that the polar decomposition of $\Log_{\pcero} \ptilde$ is $\Log_{\pcero} \ptilde=|\Log_{\pcero} \ptilde| \utilde$, where $|\Log_{\pcero} \ptilde|=\frac{\arcsin 2|[\pcero,\ptilde]|}2$ and where the partial isometry $\utilde$ is $\utilde=-\tilde{\rho_0} \vtilde^*$.
	
	In order to explore the positive part of $\Log_{\pcero}\ptilde$ we first describe $|[\pcero,\ptilde]|$ as follows. First take $\xb=\begin{psmallmatrix}	x_1\\x_2	\end{psmallmatrix}\in\kcal$ such that $\hgot(\xb)=\xb\xb^*=\begin{psmallmatrix}	x_1 x_1^*&x_1 x_2^*\\x_2 x_1^*&x_2 x_2^*	\end{psmallmatrix}=\ptilde$ and $x_1$ is positive invertible.
	 Such a choice is unique.
	Recall that the equality $|x_1|^2+|x_2|^2=1$ implies that $|x_1|=x_1$ commutes with $|x_2|$.
	
	Then if	$x_2=w|x_2|$ is the polar decomposition of $x_2$ we have 
	$
	|[\pcero,\ptilde]|^2=\begin{psmallmatrix} x_1 x_2^*x_2 x_1 &0\\ 0& x_2x_1^2x_2^*	\end{psmallmatrix}
	$
	so 
	$$
	|[\pcero,\ptilde]|=\begin{psmallmatrix} x_1 |x_2|&0\\0&wx_1|x_2| w^*	\end{psmallmatrix}.
	$$
 So we have the following expression for $|\Log_{\pcero}\ptilde|$
 $$|\Log_{\pcero}\ptilde|= \begin{psmallmatrix} \frac{\arcsin (2 x_1 |x_2|)}2&0\\0&\frac{\arcsin (2wx_1|x_2| w^*)}2	\end{psmallmatrix}.
	$$
Now write $x_1=\cos\varphi$ for a unique positive element $\varphi\in\a$ ($0\leq \varphi\leq \pi/4$, see Lemma \ref{lema norma corchete menor que 1 medio entonces p11 mayor o igual que 1 medio}), and therefore $|x_2|=\sin\varphi$. So we have proved the following result.

\begin{teo} \label{teo del angulo fi entre p0 y ptilde}
	Let $\ptilde\in\rcal$ such that  $\|[\pcero,\ptilde]\|<1/2$. Then there exists a unique element $\varphi\in \a$ ($0\leq \varphi\leq \pi/4$) such that 
	$$|\Log_{\pcero}\ptilde|= \begin{psmallmatrix}  \varphi&0\\0&w \varphi w^*	\end{psmallmatrix}.
	$$
	where $\xb=\begin{psmallmatrix}
		x_1\\x_2
	\end{psmallmatrix}$ is the element in $\kcal$ that projects on $\ptilde$ with $x_1$ positive and invertible. Here $w$ is the partial isometry of the polar decomposition $x_2=w|x_2|$.
	
	We call the positive operator $\varphi\in\a$  the {\rm{angle}} between $\pcero$ and $\ptilde$.
\end{teo} 
\begin{rem}
	Since $\Log_{\pcero}\ptilde$ directs the geodesic in $\rcal$ from $\pcero$ to $\ptilde$ in $\rcal$, then its norm is the Finsler distance from $\pcero$ to $\ptilde$, and therefore this distance is $\|\varphi\|$.
\end{rem}		
\begin{rem}
	The geometric interpretation of the commutator $[\pcero,\ptilde]$ is given by the constructions of ``projection and section'' illustrated by Figure \ref{razon doble} in exactly the same way. The correspondence $x\mapsto y$ in the picture is the classical cross ratio as defined by Zelikin in \cite{zelikin}. In our case the commutator $[\pcero,\ptilde]$ is given geometrically by the correspondence $x\mapsto w$ from $l_3$ to $l_4$ (and similarly $l_4\to l_3$). We call this correspondence the \textit{complementary cross ratio}.
\end{rem}

	\section{The Hopf fibration}\label{The Hopf fibration}
	
	In this section we will describe more properties of the Hopf fibration defined in \ref{def fibracion Hopf}. Recall that the total space is the sphere $\kcal$ in $\a^2$, the base is the Riemann sphere $\rcal$ of $\a$, the group is the unitary group $\ucal$ of the algebra $\a$ and the projection is 
	$$
	\begin{tikzcd}
		\kcal\arrow[d,"\hgot"'] \\
		\rcal 
	\end{tikzcd} \
	\ \text{ given by }	\hgot(\xb)=\xb \xb^*=\px
	$$
	
	Note that $\ucal$ operates on the right in $\kcal$ by $\xb u=\begin{psmallmatrix}
		x_1 u\\x_2 u
	\end{psmallmatrix}$ for $\xb\in \kcal$ and $u\in\ucal$.
	
	\begin{teo}\label{teo fibr hopf estruct analitica}
		The Hopf fibration $\hgot$ is a C$^\infty$ principal bundle with structure group $\ucal$.		
	\end{teo}
	\begin{proof}
		We have to prove that 
		\begin{enumerate}
			\item $\hgot$ is an C$^\infty$ map onto $\rcal$,
			\item the fibers of $\hgot$ are the $\ucal$ orbits of the action, and
			\item the map $\hgot$ has C$^\infty$ local cross sections.
		\end{enumerate}		
		1. We use the atlas defined on $\kcal$ (see Section \ref{sec cartas en K}). It will suffice to prove this only for the case of the local identification associated with $\kcal_0$ (because local identifications are obtained by just acting with $\ucal_2$ on ``basic identification'' associated to $\kcal_0$). In this identification the map $\hgot$ reads $(a,u)\mapsto \tilde p$ where $a\in\a$, $u\in\ucal$, $\tilde p=\xb\xb^*$ and $\xb=
		\begin{psmallmatrix}
			(1+a^*a)^{-1/2} u 
			\\
			a (1+a^*a)^{-1/2} u
		\end{psmallmatrix}
		$
		which is obviously smooth.
		
		The map is clearly surjective since given $\tilde p\in\rcal$, we have $\tilde p=
		\tilde u
		\begin{psmallmatrix}
			1&0\\0&0
		\end{psmallmatrix} \tilde u^*$
		and therefore
		$\tilde p=
		\tilde u
		\begin{psmallmatrix}
			1&0\\0&0
		\end{psmallmatrix} \tilde u^*=\tilde u \eb_1 \eb_1^* \tilde u^*$ so that $ \tilde p=\hgot(\tilde u\eb_1)$.
		
		2. 
		Here we will show that every fiber $\hgot^{-1}(\tilde p_\xb)$ can be identified with the unitary group $\ucal$ of $\a$ and that every projector $\tilde p_\xb\in\rcal$ is the image by $\hgot$ of one of such fibers.
		
		Let us consider first the fiber over $\pcero$ which is  $\{\tilde u\eb_1\in \kcal: \tilde u\in \ucal_2 \text{ and } \tilde u\pcero\tilde u^*=\pcero\}$. Here the equation $\tilde u\pcero\tilde u^*=\pcero$ for $\tilde u=\begin{psmallmatrix}
			u_{1,1}&u_{1,2}\\ u_{2,1}&u_{2,2}
		\end{psmallmatrix}$ is equivalent to $\begin{psmallmatrix} u_{1,1}\\ u_{2,1}\end{psmallmatrix} \begin{psmallmatrix} u_{1,1}^*& u_{2,1}^*\end{psmallmatrix} =\pcero=\begin{psmallmatrix} 1& 0\\0&0\end{psmallmatrix}$. This equation implies that $u_{1,1}\in \ucal$ and that $v_{2,1}=0$. Moreover, using that $\tilde u\in\ucal$ we can conclude that $u_{1,2}=0$ and $u_{2,2}\in \ucal$ also. Hence in this case the fiber is
		\begin{equation}
			\label{eq fla de fibra en pcero}
			\hgot^{-1}(\pcero)=\left\{\begin{psmallmatrix}
				u_{1,1}&0\\ 0&u_{2,2}
			\end{psmallmatrix}\eb_1 : u_{1,1} , u_{2,2}\in \ucal\right\}=\left\{\begin{psmallmatrix}
				u\\ 0
			\end{psmallmatrix} : u\in \ucal\right\}.
		\end{equation}
		which can be identified with $\ucal$.
		
		Now consider the general case of $\px\in \rcal$ where $\xb=\tilde v\eb_1$ with $\tilde v\in \ucal_2$ and suppose that $\hgot(\tilde w\eb_1)=\px$ with $\tilde w\in \ucal_2$ is any other element of the fiber. Then we have that $\tilde w\eb_1 \eb_1^*\tilde w^*=\px=\tilde v\eb_1 \eb_1^*\tilde v^*$ for $\tilde w\in \ucal_2$.
		Now, using that $\pcero=\eb_1 \eb_1^*=\tilde v^*\tilde w\eb_1 \eb_1^*\tilde w^*\tilde v=\tilde v^*\tilde w\pcero\tilde w^*\tilde v$ and that $\tilde v^*\tilde w\in\ucal_2$, the description \eqref{eq fla de fibra en pcero}
		proves that $\tilde v^*\tilde w= \begin{psmallmatrix}
			u_{1,1}&0\\ 0&u_{2,2}
		\end{psmallmatrix}$ for $u_{1,1}, u_{2,2}\in\ucal$. Then 
		$$
		\tilde w\eb_1=\tilde v\begin{psmallmatrix}
			u_{1,1}&0\\ 0&u_{2,2}
		\end{psmallmatrix}\eb_1=\begin{pmatrix}
			v_{1,1}u_{1,1}\\ v_{1,2} u_{1,1}
		\end{pmatrix}=\begin{pmatrix}
			v_{1,1}\\ v_{1,2} 
		\end{pmatrix}u_{1,1}\ \text{ with } u_{1,1}\in\ucal
		$$
		which proves that for every $\xb\in \kcal$ the fiber of $\px$ can be identified with $\ucal$.

		3. Consider the set of $\rcal$ given by $V_0=\hgot(\kcal_0)=\{ \xb\xb^*:\xb\in\kcal \text{ and } x_1\in \mathcal G\}\subset \rcal$, where $\kcal_0$ is the domain of the map $\psi_0$ defined in \eqref{eq def de K0 dominio carta} which is also the range of its inverse $\Psi_0:\a\otimes \ucal\to \kcal_0$ (see \eqref{eq Psi0 inversa de psi0}).
		$V_0$ is an open neighborhood of $\pcero$ where we can define a section $\sigma$ as
		\begin{equation}
			\begin{split}\label{eq seccion local de hgot}
				\sigma&:V_0\to \kcal_0 \\
				\sigma(\px)&=\Psi_0(a,1)
				=\begin{pmatrix}
					(1+a^*a)^{-1/2}  
					\\
					a (1+a^*a)^{-1/2}  
				\end{pmatrix}=\begin{pmatrix}
					1\\a   \end{pmatrix}(1+a^*a)^{-1/2}
			\end{split}
		\end{equation}
		
		for $\xb=\Psi_0(a,u)$, with $u\in\ucal$ (or equivalently for $\xb=\begin{psmallmatrix}
			x_1\\x_2\end{psmallmatrix}$ such that $x_2 x_1^{-1}=a$).
		
		Let us see first that $\sigma$ is well defined. Suppose that $\px=\xb\xb^*=\zb\zb^*=\pz$.Then $\xb=\zb u$ for $u\in \ucal$ (see Proposition \ref{prop equivs [x] = [z]}). Hence, if $x_1=r\, v$ for $r>0$ and $v\in \ucal$, then $z_1=r\,v\, u$ for $v\in\ucal$ and $z_2=x_2\,u$.
		Therefore $z_2z_1^{-1}=x_2uu^*x_1^{-1}=x_2x_1^{-1}$ which implies that $\sigma(\px)=\sigma(\pz)$.
		
		Moreover, if we compose $\sigma$ with the map $\psi_0:\kcal_0\to \a\times\ucal$ we obtain $\psi_0(\sigma(\px))=\psi_0(\Psi_0(x_2x_1^{-1},1))=(x_2x_1^{-1},1)$ which is clearly C$^\infty$ since $(\px)_{2,1} \left(\px\right)_{1,1}^{-1}=x_2 x_1^* (x_1 x_1^*)^{-1}= x_2 x_1^* (x_1^*)^{-1} x_1^{-1}=x_2 x_1^{-1}$ is an analytic function of two of the entries of $\px\in M_2(\a)$.
	\end{proof}
	
	\begin{defi}
		 Given $\ptilde\in\rcal$ we will say that $(x_1,x_2)\in\a\times \a$ is a {\rm{pair of homogeneous coordinates}} for $\ptilde$ if $\xb=\begin{psmallmatrix}
		 	x_1\\x_2
		 \end{psmallmatrix}\in \im(\ptilde)$ and there exists an invertible element $\lambda\in\a$ such that $\xb \lambda\in\kcal$. 
		\end{defi}
		 Observe that every $\ptilde\in \rcal$ has a pair of homogeneous coordinates. Also note that if $(x_1, x_2)$ and $(x'_1,x'_2)$ are pairs of homogeneous coordinates of $\ptilde$ there exists an invertible element $\lambda$ in $\a$ such that $x'_1=x_1\lambda$ and $x'_2=x_2\lambda$.
		 
		 \begin{rem}
		 	The open set $\vcal_0$ defined in \eqref{eq def V0 para p0} consists of all $\ptilde$ that have homogeneous coordinates $(x_1,x_2)$ with $x_1$ invertible.
		 \end{rem}
		 
		 \subsection{Relation between geodesic and homogeneous coordinates in $\rcal$}\label{Relation between normal and homogeneous coordinates in R}
		
		We now give an explicit expression for the relation between homogeneous coordinates and geodesic coordinates of an element $\ptilde\in U_0$.

		\begin{teo}
			Let $\ptilde\in U_0$. Consider the following diagram 
			\[
			\begin{tikzcd}
				&	 \kcal_0\subset \kcal   \arrow[rd, "\psi_0"] \arrow[d, "\hgot "] &\\
				V_0\subset	(T\rcal)_{\pcero} \arrow[r, "\Exp_\pcero"] & U_0\subset \rcal    \arrow[r, "\varphi_0"] &  \a\\
				X=\begin{psmallmatrix}	0 &a^* \\ a&0\end{psmallmatrix} \arrow[u, phantom, sloped, "\in"] \arrow[rr,mapsto] &&v\tan|a|\arrow[u, phantom, sloped, "\in"]
			\end{tikzcd}
			\] 
		\end{teo}
		where $\Log_\pcero(\ptilde)=X$ and $v$ comes from the polar decomposition $a =v|a |$. Then 
		$$
		\varphi_0(\ptilde)=v \tan(|a|).
		$$
		
		\begin{proof}
				
		Suppose that $\Exp_\pcero(X)=\ptilde$ with $X=\begin{psmallmatrix}
			0&a^* \\ a&0
		\end{psmallmatrix}$ and $\|X\|=\|a\|=\|a^*\|<\pi/4$. We know from \eqref{ec fla geodesicas en teo}, considering that in Theorem \ref{teo geods desde P(Gr(0)) con CI} we used $X=\begin{psmallmatrix}
		0&a \\ a^*&0
		\end{psmallmatrix}$, that $\ptilde=
		\begin{psmallmatrix} \cos|a|\\ a \sinc|a|	\end{psmallmatrix}\begin{psmallmatrix} \cos|a| &   (\sinc|a| a^*)	\end{psmallmatrix} $ is a possible expression of $\ptilde$ in terms of $a$.
			Now consider $\xb=\begin{psmallmatrix} x_1\\ x_2	\end{psmallmatrix} =\begin{psmallmatrix} \cos| a |\\ a (\sinc| a |)	\end{psmallmatrix}=\begin{psmallmatrix} \cos| a |\\ v (\sin| a |)	\end{psmallmatrix}$ where $v$ is the isometry of the polar decomposition $a =v|a |$. Then $x_1=\cos|a |$ is invertible since $\|a \|<\pi/4$ and $\xb^*\xb=(\cos|a |)^2+|a |^2(\sinc|a |)^2=(\cos|a |)^2+(\sin|a |)^2=1$. We can also find a unitary $\tilde u\in \ucal_2$ as in \eqref{eq: unitario de U2 para suryectividad} such that $\xb=\utilde \begin{psmallmatrix}
				1\\0
			\end{psmallmatrix}$ and then $\xb\in\kcal$.
		Therefore, $x_2x_1^{-1}=v\sin|a |(\cos|a |)^{-1}=v\tan|a |$ and hence we obtained the formula $\varphi_0(\ptilde)=v\tan|a|$.
						\end{proof}
	\begin{rem}\label{rem unoriented angle}
		 Consider the classical picture from Figure \ref{fig angulo no orientado}.		 
		\begin{figure}[h!]
			\centerline{\includegraphics[width=0.45\textwidth]{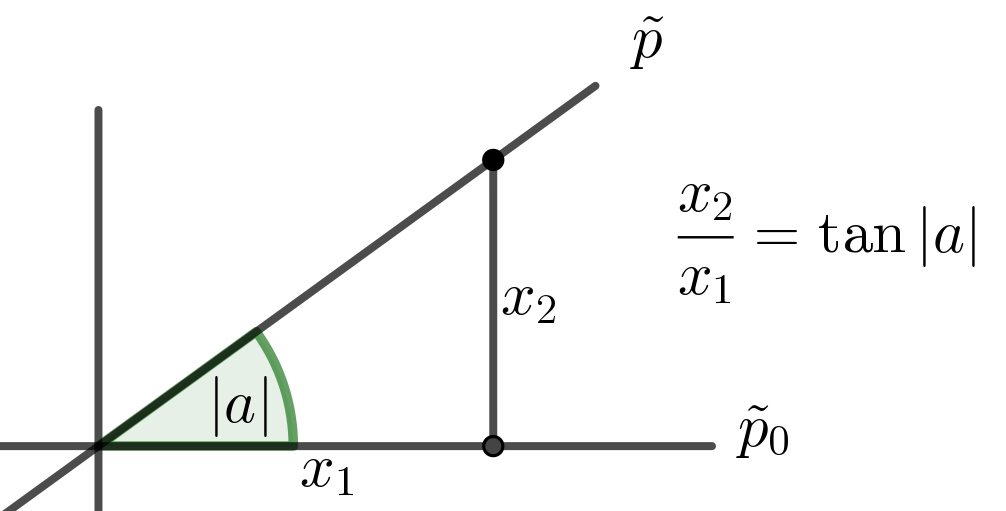}}
			\caption[short]{Unoriented angle.}
			\label{fig angulo no orientado}
		\end{figure}
		  In it we have schematically represented homogeneous coordinates $(x_1,x_2)$ for $\ptilde$. The ``affine'' coordinate $\varphi_0(\ptilde)$ is $x_2x_1^{-1}$ while the element $v\tan|a|$ is related to the geodesic coordinate of $\ptilde$ which is $X$. This suggests naming $|a|$ as the \textit{unoriented angle} between $\ptilde$ and $\pcero$ and the partial isometry $v$ becomes a ``phase'' related to the pair $(\pcero, \ptilde)$.
		  
		  Note that this angle $|a|$ coincides with the one denoted with $\varphi$ in Theorem \ref{teo del angulo fi entre p0 y ptilde}.
	\end{rem}

	\subsection{The canonical connection on the Hopf fibration}\label{The canonical connection on the Hopf fibration}
	
	We will define a C$^\infty$ \textit{horizontal distribution} $H_\xb$ of subspaces of the tangent spaces $(T\kcal)_\xb$ for $\xb\in \kcal$.
	This distribution will turn out to be invariant under the right action of $\ucal$ on $\kcal$ and consequently will define a connection on the principal bundle $\kcal\to \rcal$. We will call this connection the \textit{canonical connection} on the Hopf fibration.
	
	Observe that given $\xb\in \kcal$ the tangent space $(T\kcal)_\xb$ is described as follows
	$$
	(T\kcal)_\xb=\{\xi\in\a^2:\langle\xb,\xi\rangle \text{ is anti-self-adjoint}\}.
	$$
	Observe that every $\xb\in\kcal$ verify $\langle \xb^*,\xb\rangle=1$.

	At each point $\xb\in \kcal$ we have the vertical tangent space $V_\xb\subset (T\kcal)_\xb$ defined by 
	$$
	V_\xb=\ker(T\hgot)_\xb
	$$
	where $T\hgot$ is the tangent map. Clearly $V_\xb$ is the image of the Lie algebra of $\ucal$ under the derivative at $u=1$ of the action $u\mapsto \xb u$. Therefore $V_\xb=\{\xb a: a\in \a\text{ antiselfadjoint}\}$.
	
	Next we define the \textit{horizontal space $H_\xb \subset (T\kcal)_\xb$ at} $\xb$   for $\xb\in\kcal$ as follows
	$$
	H_\xb=\ker(\px).
	$$
	where the vectors in $H_\xb$ are considered as tangent vectors to $\a^2$ at $\xb$ (note that if $\xi\in \ker(\px)$ then $\langle\xb,\xi\rangle=0$). 
	
	Observe that $(T\kcal)_\xb=V_\xb\oplus H_\xb$. It is also clear that 
	the map $TR_u:(T\kcal)_\xb\to (T\kcal)_{\xb u}$ (where $R_u$ is the right multiplication and where $TR_u$ is the tangent map of $R_u$) satisfies
	$$
	(TR_u)_\xb (H_\xb)=H_{\xb u}.
	$$
	This completes the statement at the beginning of this paragraph about the definition of the canonical connection on the Hopf fibration.
	
	\begin{rem}
		Clearly, the (left) action of $\ucal_2$ on $\kcal$ preserves the decomposition $(T\kcal)_\xb=V_\xb\oplus H_\xb$.
	\end{rem}	
	We finish this section describing the tangent map $(T\hgot)_\xb:(T\kcal)_\xb\to(T\rcal)_\px$ for $\xb \in\kcal$. Given $\xi\in (T\kcal)_\xb$ we clearly have that $T\hgot (\xi)=X=\xi \xb^*+\xb\xi^*$. Also note the identity $X\px=(1-\px) X$.
	
	\subsection{The structure morphism $\kappa: \mathfrak{R}\to T\kcal$}\label{subsecc de kapa}
	
	Define the vector bundle $\mathfrak{R}\to \kcal$ as the induced vector bundle 
	\[
	\begin{tikzcd}
		\hgot^* (T\rcal) \arrow{r}  \arrow{d}{} &T\rcal \arrow{d}{} \\   
		\kcal \arrow[swap]{r}[swap]{\hgot} & \rcal
	\end{tikzcd}
	\]
	where we write $\mathfrak{R}$ for $\hgot^* (T\rcal)$ as a bundle over $\kcal$. With this notation we define the \textit{structure morphism} $\kappa$ as a vector bundle morphism 
	$$
	\begin{tikzcd}
		\mathfrak{R} \arrow[r, "\kappa"] \arrow[dr]
		&  {T\kcal} \arrow[d]\\
		& \kcal
	\end{tikzcd}
	$$
	where $\kappa_\xb(X)=X\xb$ for each $\xb\in\kcal$, $X\in (T\rcal)_\px$ (notice that $X\in M_2(\a)$, $\xb\in \a^2$ and therefore $X\xb\in \a^2$). Observe that $\kappa_\xb(X)\in \ker(\px)=  H_\xb$. Also notice that $(T\hgot)_\xb(\kappa_\xb X)=X$, because of the identity $X\xb\xbstar+\xb\xbstar X=X$ (observe that $X$ is self-adjoint).
	
	We remark that the morphism $\kappa$ has the following equivariance property 
	$$
	\kappa_{\xb u}(X)=\left(\kappa_\xb(X)\right) u
	$$
	where $u\in\ucal$. This equivariance shows a way of constructing the tangent bundle $T\rcal$ out off the principal bundle $\kcal\to\rcal$ and the co-tautological bundle $\mathcal{T}'$ (see \eqref{eq T y Tprima}).
	
	The following schematic picture illustrates our constructions 
	\\
	
	\centerline{\includegraphics[width=0.75\textwidth]{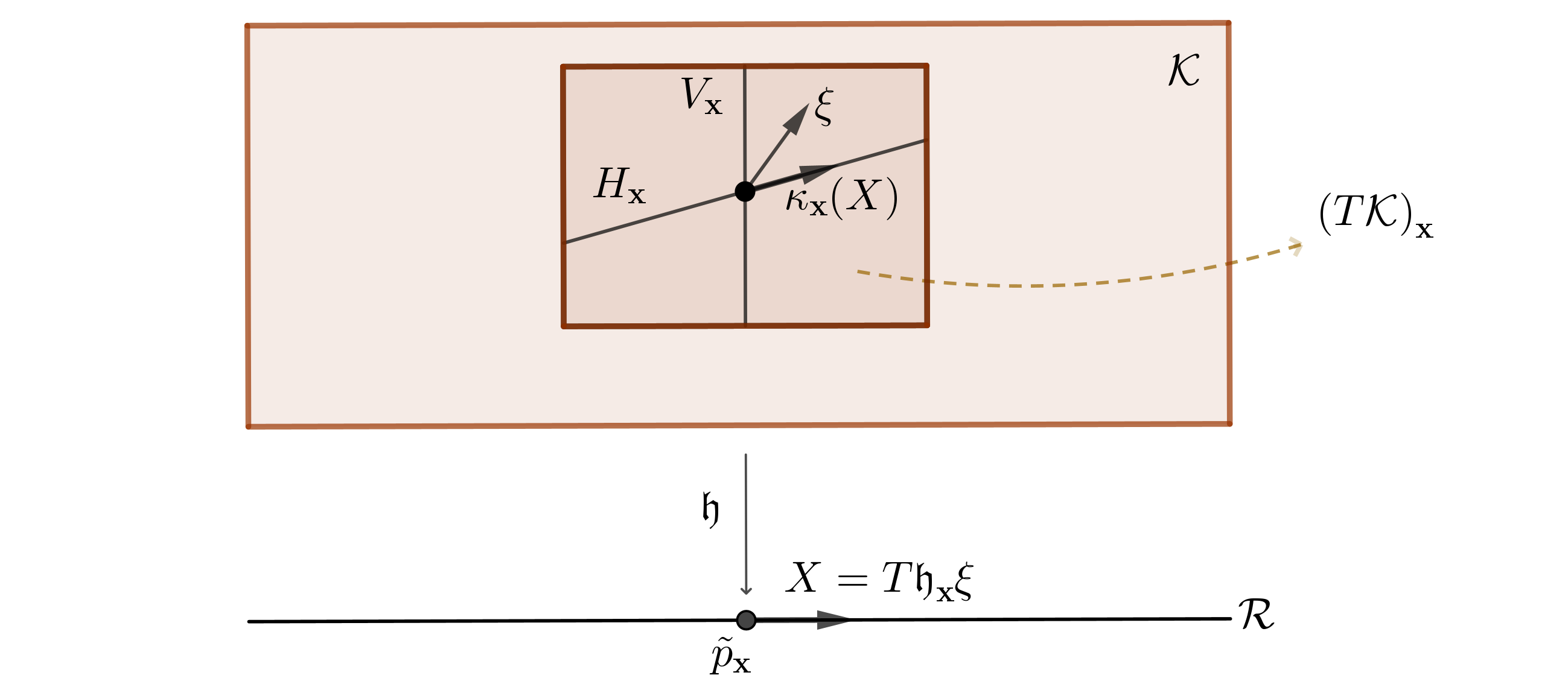}}
	
	where the inner rectangle represents the tangent space $(T\kcal)_\xb$.
	
	\subsection{The Finsler metric on $\rcal$ and the structure form $\kappa$}\label{subsec metrica de Finsler}
	Recall that $\a^2$ is a Hilbert C$^*$-module over $\a$ (acting on the right) in the usual way defining $\langle \xb,\yb\rangle=x_1^*y_1+x_2^*y_2$. Then we have the following 
	\begin{teo}
		Let $X\in(T\rcal)_{\tilde p}$, $\xb\in\kcal$, $\hgot(\xb)=\tilde p$. Then
		$$
		\|X\|=\|\kappa_\xb(X)\|=\|X\xb\|.
		$$
		Here $\|X\|$ is the Finsler norm in $\rcal$ of the tangent vector $X$ (i.e. the usual norm of the self-adjoint matrix $X\in M_2(\a)$)  whereas $\|\kappa_\xb(X)\|$ stands for the norm of $\kappa_\xb(X)$ as an element of the C$^*$ $\a$-module $\a^2$.
	\end{teo}
	\begin{proof}
		Suppose first that $\tilde p=\pcero$. In this case, since $X\in (T\rcal)_\pcero$, we have that $X=\begin{psmallmatrix}
			0&a\\a^*&0\end{psmallmatrix}$ for $a\in\a$. Then
			$$
			\|X\|_{M_2(\a)}=\|a\|=\|X \begin{psmallmatrix}
				1\\0\end{psmallmatrix}\|_{\a^2}=\|X\begin{psmallmatrix}
				1 u\\0\end{psmallmatrix}\|_{\a^2}=\|X\xb\|_{\a^2}
			$$
			where $u\in \ucal$, $\xb\in \kcal$ with $\hgot(\xb)=\xb\xb^*=\pcero$, and $\|\,\|_{\a^2}$ is the norm of the Hilbert C$^*$-module $\a^2$.
			
			The general case follows using that given $\ptilde\in\rcal$, there is $\zb\in\kcal$ such that  $\ptilde=\pz=\zb\zb^*=\utilde \begin{psmallmatrix}
				1\\0\end{psmallmatrix}$ with $\utilde\in\ucal$. And every element of $(T\rcal)_\ptilde$ is of the form $\tilde u X\utilde^*$ for $X\in (T\rcal)_\pcero$. Hence
			$\|\utilde X\utilde^*\|_{M_2(\a)}=\|X\|_{M_2(\a)}$. And for $\yb$ such that $\hgot(y)=\ptilde$ we have that  $\yb=\zb v$ with $v\in\ucal$. Then we obtain that
			\begin{equation*}
				\begin{split}
			\|\left(\utilde X\utilde^*\right) \yb\|_{\a^2}&=
			\|\left(\utilde X\utilde^*\right) \zb v
			\|_{\a^2}=
			\|\left(\utilde X\utilde^*\right) 
			\utilde \begin{psmallmatrix}
				1\\0\end{psmallmatrix} v
			\|_{\a^2}
			=	\|\utilde X\begin{psmallmatrix}
				1\\0\end{psmallmatrix} v
			\|_{\a^2}
		\\
			&
			=\|v^*\begin{psmallmatrix}
				1&0\end{psmallmatrix} X^*\utilde^*\utilde X\begin{psmallmatrix}
				1\\0\end{psmallmatrix} v
			\|^{1/2}
			=\|\begin{psmallmatrix}
				1&0\end{psmallmatrix} X^* X\begin{psmallmatrix}
				1\\0\end{psmallmatrix} 
			\|^{1/2}
			=\|X\|_{M_2(\a)}
			\end{split}
			\end{equation*}
			where the last equality was proved in the case of $\ptilde=\pcero$.
	\end{proof}
	
	\section {Examples}
	
	\subsection{The finite dimensional case}
	
	Consider $\a=M_n(\CC)$. The case where $n=1$ is the classical Riemann $\rcal$ sphere and the classical Hopf fibration $\kcal\to \rcal$ over the Riemann sphere (see \ref{def fibracion Hopf}). In this case $\rcal$ is the one dimensional complex projective line $\mathbb{P}(\CC)$ (homeomorphic to $S^2$) and $\kcal$ is the unit sphere in $\CC^2$ which is homeomorphic to $S^3$.
	
		The case $n>1$ involves the non commutative C$^*$-algebra $M_n(\CC)$ of operators on $H=\CC^n$. Here $M_2(\a)$ is naturally identified with $M_{2n}(\CC)$ operating on $H\oplus H$ which is naturally identified with $\CC^{2n}$.
		Also $\pcero$ is the orthogonal projection in $\CC^{2n}$ onto $\CC^n\subset \CC^{2n}$ as the subspace defined by $z_{n+1}=z_{n+2}=\dots=z_{2n}=0$.
		Therefore the orbit $\rcal$ of $\pcero$ by the action of $\ucal_2\subset M_2(\a)$ can be identified with the classical Grassmann manifold $\text{Grass}_{n,2n}(\CC)$ of all $n$ dimensional subspaces of $\CC^{2n}$.
		
		We now describe the sphere $\kcal$ in $\a^2$ corresponding to the present situation. We have that
		$$
		\kcal=\{\xb=\begin{psmallmatrix}
			x_1\\x_2
		\end{psmallmatrix}:\CC^n\to\CC^n\oplus \CC^n : \xb \text{ is an isometry} \}.
		$$
	Observe that $\xb^*=\begin{psmallmatrix}
		x_1^*&x_2^*
	\end{psmallmatrix}:\CC^{2n}\to\CC^n$ and $\xb\xb^*$ is an orthogonal projection in $\CC^{2n}$ so that $\hgot:\kcal\to\rcal$ is given by the usual formula. 
	The space $\kcal$ may be identified with the usual Stieffel manifold $St_{n,2n}$ of orthogonal $n$-frames in $\CC^{2n}$ and $\hgot$ is therefore identified to the usual projection $St_{n,2n}\to \text{Grass}_{n,2n}$.
		
	In this context the open set $\vcal_0$, domain of the principal chart, consists of all orthogonal projections $\tilde p\in\rcal$ such that $\im \tilde p$ is the graph of a linear map $a:\CC^n\to\CC^n$ and $\varphi_0(\tilde p)=a$.
	
	\begin{rem}
		Notice that $\vcal_0$ is dense in $\rcal$. In the standard CW-decomposition of $\rcal$, $\vcal_0$ is the \textit{top cell} and has (real) dimension $4n^2$. See for example \cite{milnor} for the real case. The complex case is similar.
	\end{rem}

	\subsection{Bounded and unbounded operators}\label{secc The case of unbounded operators}
In this subsection we present the closed operators on a Hilbert space $H$ as elements of the Riemann sphere of the algebra $\a=B(H)$.	

%
%
%
%
	
	For a densely defined closed operator $T:\text{Dom}(T)\to H$, it can be proved that its orthogonal projection $P_{\Gr(T)}$ onto the graph of $T$ belongs to the Riemann sphere of $\a=L(H)$. These statements are formalized in the following result where we also provide formulas for these orthogonal projections.
	
	\begin{prop}\label{prop los proy sobre graficos de no acotados estan en R}
		Let $H$ be  a Hilbert space and $T:\text{Dom}(T)\to H$ a densely defined operator with closed graph. Then the orthogonal projection $\pp_T=P_{\Gr(T)}$ over the graph $\Gr(T)=\{(h,T(h)):h\in \text{Dom}(T)\}\subset H\times H$ belongs to the unitary orbit of $\left(\begin{smallmatrix}1&0\\0&0\end{smallmatrix}\right)$ which is the Riemann sphere $\mathcal{R}$ of the algebra $B(H)$.
		Moreover, $P_{\Gr(T)}$ can be written as
		\begin{equation}\label{eq proy asociada a proy grafico op no acotado}
			\begin{split}
				P_{\Gr(T)}&=\begin{pmatrix}
					1\\ T
				\end{pmatrix}(1+T^*T)^{-1}\begin{pmatrix}
					1&T^*
				\end{pmatrix} \\
				&=\begin{pmatrix}1&T^*\\ T&TT^*
				\end{pmatrix}\begin{pmatrix}
					1+T^*T&0\\0&1+TT^*
				\end{pmatrix}^{-1}\\
				&= \begin{pmatrix}
					(1+T^*T)^{-1}&(1+T^*T)^{-1}T^*\\ T(1+T^*T)^{-1}&T(1+T^*T)^{-1}T^*
				\end{pmatrix}.	
				\end{split}
		\end{equation}
	Observe that all the entries of the last matrix are bounded operators.	
	\end{prop}
	\begin{proof}
		Consider the operator $\tilde T$ defined by
		\begin{equation}\label{def f tilde}
			\tilde T=\begin{pmatrix}
				0&-T^*\\ T&0
			\end{pmatrix}.
		\end{equation}
		Observe that since $T$ is a closed and densely defined operator on a Hilbert space then $T^*$ is also closed and densely defined (see \cite[Theorem 5.3]{weidmann}
		) and $T^{**}=T$ \cite[Theorem 1.8]{schmudgen}. Moreover, $\tilde T^*=-\tilde T$ and then $\text{Dom}(\tilde T)=\text{Dom}(\tilde T^*)$.
		
		Now we will consider the norms. Given $\left(\begin{smallmatrix}
			\xi\\ \eta 
		\end{smallmatrix}\right)\in\text{Dom}(\tilde T)$ we have that
		$$
		\tilde T \left(\begin{smallmatrix}
			\xi\\ \eta 
		\end{smallmatrix}\right)=-\tilde T^*\left(\begin{smallmatrix}
			\xi\\ \eta 
		\end{smallmatrix}\right)
		$$
		and therefore they have the same norm. This proves that $\tilde T$ is a normal operator in $H\times H$ (see \cite[Section 5.6]{weidmann}). Then since this implies that also $1+\tilde T$ is normal and then is invertible with a bounded inverse. This follows considering that 	$\tilde T^*=-\tilde T$ and then using the functional calculus of the self-adjoint operator $i\tilde T$.		 
		
		Now consider the polar decomposition (see 
		\cite[Theorem 7.20]{weidmann})
		$$
		1+\tilde T= US
		$$
		where $U$ is a unitary operator since $1+\tilde T$ is invertible. This follows because the range of $1+\tilde T$ is dense (see \cite[Theorem 5.42]{weidmann}) and hence $U$ is onto, and since $1+\tilde T$ is injective, then $U$ must be injective.

		Now let us analyze the operator $S$. Using the same reference cited above, $S$ can be written as  
		\begin{equation}\label{eq fla. para s}
			\begin{split}
				S&=|1+\tilde T|=\left((1+\tilde T)^*(1+\tilde T)\right)^{1/2}=\left((1+\tilde T^*)(1+\tilde T)\right)^{1/2}\\&
				=\left((1-\tilde T)(1+\tilde T)\right)^{1/2}=\left(1-\tilde T^2\right)^{1/2}.
			\end{split}
		\end{equation}
		
		And then, using \eqref{def f tilde}, a direct computation gives that 
		$S^2=|1+\tilde T|^2=1+\left(\begin{smallmatrix}
			T^*T&0\\0&T T^*
		\end{smallmatrix}\right)=\left(\begin{smallmatrix}
			1+T^*T&0\\0&1+T T^*
		\end{smallmatrix}\right)$
		and hence
		$$
		S=|1+\tilde T|=\begin{pmatrix}
			(1+T^*T)^{1/2}&0\\0& (1+TT^*)^{1/2}
		\end{pmatrix}.
		$$
		Then $S$ is invertible with bounded inverse (see \cite[Proposition 3.18]{schmudgen})
		and we can write
		$$
		U=\begin{pmatrix}
			1&-T^*\\T&1
		\end{pmatrix}\begin{pmatrix}
			(1+T^*T)^{-1/2}&0\\0& (1+TT^*)^{-1/2}
		\end{pmatrix}. 
		$$
		Then the first column of $U$ is $\left(\begin{smallmatrix}
			1\\T
		\end{smallmatrix}\right)(1+T^*T)^{-1/2}=\left(\begin{smallmatrix}
			(1+T^*T)^{-1/2}\\ T(1+T^*T)^{-1/2}
		\end{smallmatrix}\right)$ with an invertible first coordinate and second coordinate $Z_T=T(1+T^*T)^{-1/2}$, which is usually called the bounded transform of $T$. 
		
		Then it follows that $\a=(1+T^*T)^{-1/2} \a$ and, hence
		$
		\left[ \begin{pmatrix}
			1\\ T
		\end{pmatrix} (1+T^*T)^{-1/2}\right]=\left[\begin{pmatrix}
			1\\T
		\end{pmatrix}\right] =  \left( \begin{smallmatrix}
			1\\ T
		\end{smallmatrix}\right)\a  $ 
		that can be identified with $\Gr(T)$. 
		
		
		Then the corresponding orthogonal projection $\ptilde_T=P_{\Gr(t)}$ onto the graph $\Gr(T)$ belongs to $\rcal$ and can be written as
		\begin{equation*}
			\begin{split}
				\pp_T&=\begin{pmatrix}
					1\\T
				\end{pmatrix}(1+T^*T)^{-1/2}(1+T^*T)^{-1/2}\begin{pmatrix}
					1&T^*
				\end{pmatrix}
				=\begin{pmatrix}
					(1+T^*T)^{-1}&(1+T^*T)^{-1}T^*\\f(1+T^*T)^{-1}&T(1+T^*T)^{-1}T^*
				\end{pmatrix}\\
				& =\begin{pmatrix}1&f^*\\ T&TT^*
				\end{pmatrix}\begin{pmatrix}
					1+T^*T&0\\0&1+TT^*
				\end{pmatrix}^{-1}.	
			\end{split}
		\end{equation*}
		\end{proof}	
	The following facts will be useful to establish the existence of minimal geodesics between graphs of operators.
\begin{defi}
	\label{def cograph}
	The \textit{inverse graph} (see \cite{kato}) of a densely defined operator $T$ on $D(T)\subset H$ is given by
	\begin{equation} \label{ec def cograph}
		\mathop{invGr}(T)=\{(Tx,x): x\in \text{D}(T)\}.	
	\end{equation}
\end{defi}

\begin{lem}\label{lema Graf T perp}
	If $T:\mathcal{D}(T)\to H$ is a densely defined closed operator on $\dcal(T)\subset H$, then
	\begin{equation*}
		\begin{split}
			Gr(T)^\perp&=\{(-T^*x,x):x\in \mathcal{D}(T^*)\}\\
			&= invGr(-T^*).
		\end{split} 
	\end{equation*}
\end{lem}
\begin{proof}
	We can use the unitary operator $V:H\oplus H\to H\oplus H$ defined by $V(x,y)=(-y,x)$ to write $Gr(T^*)=V(Gr(T)^{\perp})$ (see \cite[Lemma 1.10]{schmudgen}). Then, since $V^2=-I$ we can write $Gr(T)^\perp=-V^2(Gr(T^\perp))=-V(Gr(T^*))$ and therefore
	\begin{equation}
		\begin{split}
			Gr(T)^\perp&=-V\left(\left\{(x,T^*x):x\in\mathcal{D}(T^*)\right\}\right)
			=-\left\{(-T^*x,x):x\in\mathcal{D}(T^*)\right\}\\
			&=\left\{(-T^*x,x):x\in\mathcal{D}(T^*)\right\}.
		\end{split}
	\end{equation}
\end{proof}

\begin{prop}
	Let $S,T$ be bounded operators acting in $H$. 
	\begin{enumerate}
		\item
		There exists a (minimal) geodesic of $\mathcal{R}$ joining $P_{\Gr(S)}$ and $P_{\Gr(T)}$ if and only if $\dim \ker(1+T^*S)=\dim \ker(1+T^*S)$. The minimal geodesic is unique if and only if these subspaces are trivial.
		\item
		If $S^*=S$ and $T^*=T$, the global unitary isomorphism $\Omega$ of $H\times H$ given by
		$$
		\Omega (h_1, h_2)=(h_2,-h_1)
		$$
		maps $\ker(1+TS)$ onto $\ker(1+ST)$. In particular,  there always exist a minimal geodesic of $\mathcal{R}$ joining $P_{\Gr(S)}$ and $P_{\Gr(T)}$. 
	\end{enumerate}
\end{prop} 
\begin{proof}
	Note that $\Gr(T)^\perp=\{(-T^*g,g): g\in  H\}$, therefore a pair $(h,Sh)\in \Gr(S)$ belongs to $\Gr(T)^\perp$ if and only if there exists $g\in H$ such that $h=-T^*g$ and $Sh=g$. Therefore $h=-T^*g=-T^*Sh$, i.e., $h\in \ker(1+T^*S)$. Conversely, if $h\in \ker(1+T^*S)$, then $(h,Sh)=(-T^*Sh,Sh)\in \Gr(T)^\perp$. Then 
	$$
	\Gr(S)\cap\Gr(T)^\perp=\{(h,Sh): h\in \ker(1+T^*S)\},
	$$
	and  $\dim \Gr(S)\cap\Gr(T)^\perp=\dim \ker(1+T^*S)$. Similarly, $\Gr(T)\cap\Gr(S)^\perp=\{(g,Tg): g\in \ker(1+S^*T)\}$ with the same dimension as $\ker(1+S^*T)$. The proof follows recalling that given two subspaces $V$ and $W$, the necessary a sufficient condition for the existence of a minimal geodesic joining the orthogonal projections $P_V$ and $P_W$ is the equality of the dimensions of $V\cap W^\perp$ and $V^\perp\cap W$; and that the minimal geodesic is unique if and only if these intersections are trivial (see \cite[Theorem 4.5]{survey proyectores}).
	
	Suppose now that $S$ and $T$ are self-adjoint. Note that  $h\in \ker(1+TS)$, means that $\Omega(h,Sh)=(Sh,-h)=(Sh,TSh)$ belongs to $\Gr(T)$, with $Sh\in \ker(1+ST)$:
	$STSh=S(TSh)=S(-h)=-Sh$. That is $\Omega$ maps $\Gr(S)\cap\Gr(T)^\perp$ into $\Gr(T)\cap\Gr(S)^\perp$. Similarly, $\Omega$ maps $\Gr(T)\cap\Gr(S)^\perp$ into $\Gr(S)\cap\Gr(T)^\perp$. Note that $\Omega^2=-1$.
\end{proof}
If $S$ or $T$ are non self-adjoint, there may not exist geodesics joining their graphs, consider the following example:
\begin{ejem}
	Consider ${\bf S}$ the (unilateral) shift operator in $\ell^2$: ${\bf S}(x_1,x_2, \dots)=(0,x_1, x_2, \dots)$. Put $S=-2{\bf S}$ and $T=1$. Then $\ker(1+T^*S)=\ker(1-2{\bf S})=\ker(\frac12-{\bf S})=\{0\}$ (the shift has no eigenvalues). On the other hand $\ker(1+S^*T)=\ker(1-2{\bf S}^*)=\ker(\frac12-{\bf S}^*)$ which has dimension $1$. Therefore $\Gr(1)=\{(x,x): x\in\ell^2\}$ and $\Gr(-2{\bf S})=\{(y,-2{\bf S}y): y\in\ell^2\}$ cannot be joined by a geodesic of $\mathcal{R}$.
\end{ejem}
%
\subsection{The unique minimal geodesic from $\pcero$ to the graph of a closed operator}

Let us describe explicitly the minimal geodesic $\gamma$ of $\mathcal{R}$ with $\gamma(0)=p_0=P_{\Gr(0)}$ and $P_{\Gr(T)}$, for $f:D(T)\subset H\to H$ a closed operator. Recall from (\ref{eq proy asociada a proy grafico op no acotado}) the formula of the projection $P_{\Gr(T)}$:
$$
P_{\Gr(T)}=\left( \begin{array}{cc} (1+T^*T)^{-1} & (1+T^*T)^{-1}T^* \\ T(1+T^*T)^{-1} & T(1+T^*T)^{-1}T^*\end{array}\right).
$$
Let $T=V|T|$ be the polar decomposition of $T$, where $|T|$ is (a possibly unbounded) non- negative self-adjoint operator, and $V: \overline{R(T^*)}\to\overline{R(T)}$ is a partial isometry.

\begin{teo}
	With the current notations, we have that
	\begin{equation}\label{log grafica de f}
		\gamma(t)=e^{itZ}p_0e^{-itZ}, \hbox{ for } Z=\left(\begin{array}{cc} 0 & i \arctan(|T|)V^* \\ -i V \arctan(|T|) & 0 \end{array}\right).
	\end{equation}
\end{teo}
\begin{proof}
	To verify (\ref{log grafica de f}), let us compute the even and odd powers of $itZ$.
	Note that 
	$$
	(iZ)^2=\left(\begin{array}{cc} - \arctan(|T|) V^*V\arctan(|T|) & 0 \\ 0 & - V(\arctan(|T|))^2V^*\end{array}\right).
	$$
	Since $V$ is a partial isometry with initial space $\overline{R(|T|)}$ and final space $\overline{R(T)}$, and $\arctan(T)$ is a continuous  function with $\arctan(0)=0$,  it follows that $V^*V=P_{\overline{R(|T|)}}$, and thus $V^* V \arctan(|T|)=\arctan(|T|)V^*V=\arctan(|T|)$. Therefore we have 
	$$
	(iZ)^2=\left(\begin{array}{cc} - \left(\arctan(|T|)\right)^2  & 0 \\ 0 & - V\left(\arctan(|T|)\right)^2V^*\end{array}\right).
	$$
	Similarly,
	$$
	(iZ)^{2k}=(-1)^k \left(\begin{array}{cc} - \left(\arctan(|T|)\right)^{2k}  & 0 \\ 0 & - V\left(\arctan(|T|)\right)^{2k}V^*\end{array}\right) .
	$$
	The odd powers of $iZ$: $(iZ)^3$ equal
	$$
	\left(\begin{array}{cc} - \left(\arctan(|T|)\right)^2  & 0 \\ 0 & - v\left(\arctan(|T|)\right)^2V^*\end{array}\right)\left(\begin{array}{cc} 0 & - \arctan(|T|) V^*   \\  V\arctan(|T|) & 0\end{array}\right)
	$$
	$$
	=\left(\begin{array}{cc} 0 & \arctan(|T|) V^*V \left( \arctan(|T|)\right)^2 V^* \\ -V\left(\arctan(|T|)\right)^3 & 0 \end{array}\right)
	$$
	$$
	=\left(\begin{array}{cc} 0 & \left( \arctan(|T|)\right)^3 V^* \\ -V\left(\arctan(|T|)\right)^3 & 0 \end{array}\right).
	$$
	Similarly,
	$$
	(iZ)^{2k+1}=(-1)^k \left(\begin{array}{cc} 0 & -\left(\arctan(|T|)\right)^{2k+1}V^* \\ V\left(\arctan(|T|)\right)^{2k+1} & 0 \end{array}\right).
	$$
	Therefore
	$$
	e^{iZ}=\left( \begin{array}{cc} \cos\left(\arctan(|T|)\right) & -\sin\left(\arctan(|T|)\right)
		\\ \sin\left(\arctan(|T|)\right) & \cos\left(\arctan(|T|)\right)\end{array}\right).
	$$
	Notice the functional identities $\cos(\arctan(t))=\frac{1}{\sqrt{1+t^2}}$ and $\sin(\arctan(t))=\frac{t}{\sqrt{1+t^2}}$. 
	Then (using that $T=V|T|$ and $T^*=|T|V^*$), $e^{iZ}$ equals
	$$
	\left( \begin{array}{cc} (1+|T|^2)^{-1/2} & -(1+|T|^2)^{-1/2}|T|V^* \\ V|T| (1+|T|^2)^{-1/2} & V(1+|T|^2)^{-1/2}V^* \end{array}\right)=
	\left( \begin{array}{cc} (1+|T|^2)^{-1/2} & -(1+|T|^2)^{-1/2}f^* \\ T (1+|T|^2)^{-1/2} & V(1+|T|^2)^{-1/2}T^* \end{array}\right).
	$$
	Then, after straightforward computations, $e^{iZ}p_0 e^{-iZ}$ equals
	$$
	\left( \begin{array}{cc} (1+|T|^2)^{-1/2} & -(1+|T|^2)^{-1/2}T^* \\ T (1+|T|^2)^{-1/2} & V(1+|T|^2)^{-1/2}V^* \end{array}\right)\left(\begin{array}{cc} 1 & 0 \\ 0 &  0 \end{array} \right) \left( \begin{array}{cc} (1+|T|^2)^{-1/2} & (1+|T|^2)^{-1/2}T^* \\ -T (1+|T|^2)^{-1/2} & T(1+|T|^2)^{-1/2}V^* \end{array}\right)
	$$
	$$
	=\left( \begin{array}{cc} (1+|T|^2)^{-1/2} & (1+|T|^2)^{-1/2}T^* \\ T (1+|T|^2)^{-1/2} & f(1+|T|^2)^{-1/2}T^* \end{array}\right)=P_{\Gr(T)},
	$$
	as claimed.
\end{proof}

Note that if $T$ is bounded, then $\|Z\|=\|\arctan(|T|)\|=\arctan(\|T\|)<\pi/2$, while if $T$ is unbounded, $\|Z\|=\|\arctan(|T|)\|=\pi/2$.

	\subsection{Bounded deformations of unbounded operators}\label{bounded deformations}
	
	In this section we consider operators $T$ on a Hilbert space $H$.
	
	\begin{defi}\label{def bounded deformation}
		A {\rm{bounded deformation}} of an unbounded closed operator $T$ on $H$ is a family $\{T_t\}_{t\in[0,\alpha)}$, with $\alpha >0$ of bounded operators $T_t$ such that
		\begin{itemize}
			\item $t\mapsto T_t$ is continuous in the norm topology
			
			\item $\lim_{t\to\alpha^{-}} \ptilde_t=\ptilde_T$ where $\ptilde_t$, $\ptilde_T$ are in the Riemann sphere $\rcal$ of the algebra $B(H)$, $\ptilde_t$ is the orthogonal projection on $\Gr(T_t)$, $\ptilde_T$ is the orthogonal projection on $\Gr(T)$  (cf Proposition \ref{prop los proy sobre graficos de no acotados estan en R}) and where the limit is taken in the Finsler metric of the Riemann sphere $\rcal$.
		\end{itemize}	

		In particular if the bounded deformation  $\{T_t\}_{t\in[0,\alpha)}$ of the unbounded operator $T$ satisfies the condition
		$$
		\dist(\ptilde_{t_0},\ptilde_\alpha)=\length\ \ptilde_t|_{t_0}^\alpha, \text{ for every } t_0\in [0,\alpha)
		$$
		we will call it an {\rm{optimal bounded deformation}}. Here $	\dist(\ptilde_{t_0},\ptilde_\alpha)$ stands for the Finsler distance in $\rcal$ and $\length\ \ptilde_t|_{t_0}^\alpha$ means the Finsler length of the curve where we write $\ptilde_\alpha$ for $\ptilde_T$.
	\end{defi}
 	In Theorem \ref{teo las geods de PGr(0) a PGr(T) son unicas y son proy sobre graficos} and Corollary \ref{coro deformation length dist} we construct a specific optimal bounded deformation of any unbounded operator $T$ on $H$.
	
 	\begin{rem}\label{rem PinvGr(T) esta en R}
	Observe that for an operator $T$, $P_{\invGr(T)}=\begin{psmallmatrix}	0&1\\1&0 \end{psmallmatrix} P_{\Gr(T)} \begin{psmallmatrix}	0&1\\1&0 \end{psmallmatrix}$ holds, which implies that $P_{\invGr(T)}\in \rcal$.
	\end{rem}
	
	\begin{teo}\label{teo las geods de PGr(0) a PGr(T) son unicas y son proy sobre graficos}
		Let $H$ be a Hilbert space, $Gr(0)=H\oplus\{0\}$ the graph of the null operator and $Gr(T)$ the graph of a densely defined closed operator $T$ with domain $\mathcal{D}(T)$.
		
		The unique minimal geodesic $\gamma:[0,1]\to \Grass(H\oplus H)$ such that $\gamma(0)=P_{\text{Gr}(0)}$ and $\gamma(1)=P_{\text{Gr}(T)}$  consists of orthogonal projections onto the graphs 
		$$
		\gamma(t)=P_{\text{Gr}(A(t))}, \text{ with } A(t)= t a^* (\sinc|t a^*|) (\cos|t a^*|)^{-1}=v \tan|t a^*|\in B(H),
		$$ 
		for $t\in[0,1)$ and $v$ the partial isometry of the polar decomposition of $a^*=v|a^*|$, with $\|a\|\leq\pi/2$.
	\end{teo}
	
	\begin{proof} Note that $\text{ran}(P_{\text{Gr}(0)})=\text{Gr}(0)=H\oplus \{0\}$,  $\ker(P_{\text{Gr}(0)})=\text{Gr}(0)^{\perp}$, $\text{ran}(P_{\text{Gr}(T)})=\text{Gr}(T)=\{(x,Tx):x\in\text{Dom}(T)\}$ and $\text{ran}(P_{\text{Gr}(T)})=Gr(T)^\perp$.
		Observe that $\text{Gr}(0)^\perp\cap \text{Gr}(T)=\left(\{0\}\oplus H\right)\cap \{(x,Tx):x\in\text{Dom}(T)\}=\{(0,0)\}$ for any $T$.
		Then we only need to prove that $\text{Gr}(0)\cap \text{Gr}(T)^\perp= \{(0,0)\}$ and use \cite[Theorem 4.5]{survey proyectores}.
		With this objective, using Lemma \ref{lema Graf T perp}, we express $\text{Gr}(T)^\perp=\{(-T^*x,x):x\in \mathcal{D}(T^*)\}$ and then obtain $\text{Gr}(0)\cap \text{Gr}(T)^\perp= \left(H\oplus\{0\}\right)\cap \{(-T^*x,x):x\in \mathcal{D}(T^*)\}=\{(0,0)\}$ which proves the uniqueness.

		Let $\gamma:[0,1]\to \Grass(H\oplus H)$ be the unique geodesic that joins $P_{\text{Gr}(0)}=\begin{psmallmatrix}1&0\\0&0\end{psmallmatrix}$ with $P_{\text{Gr}(T)}$ such that $\dot\gamma(0)=\begin{psmallmatrix}0&a\\a^*&0\end{psmallmatrix}$ (see \cite[Proposition 2.9]{survey proyectores}). This $\gamma$ is of the form (see \eqref{ec fla geodesicas en teo} and Theorem \ref{teo geods desde P(Gr(0)) con CI}).
		$$\gamma(t)=
		\begin{psmallmatrix} \cos|t a^*|\\ (\sinc|t a|)ta^*	\end{psmallmatrix}\begin{psmallmatrix} \cos|t a^*| & ta(\sinc|t a|)	\end{psmallmatrix}
		$$
		for $t\in[0,1]$.
		Moreover, the geodesics can be chosen to satisfy that $\|\begin{psmallmatrix}0&a\\a^*&0\end{psmallmatrix}\|=\|a\|\leq \pi/2$  (see  \cite[Proposition 3.1, Theorem 3.2]{survey proyectores}). 
		For all $t\in[0,1)$, the vectors $\xb(t)=\begin{psmallmatrix} \cos|t a^*|\\ (\sinc|t a|)ta^*	\end{psmallmatrix}\in\kcal_0$ because $\tilde v(t)=\begin{psmallmatrix} \cos|t a^*|& -\sinc|ta^*|ta \\ \sinc|ta|ta^*&\cos|ta|\end{psmallmatrix}\in\ucal_2$ and $\tilde v(t) \begin{psmallmatrix}
			1\\0	\end{psmallmatrix}=\xb(t)$
			and satisfy that $\xb(t)\xb(t)^*\in \vcal_0$ (see \eqref{eq def V0 para p0}) since $\cos|ta^*|$ is invertible if $\|a\|\leq\pi/2$.
		
		Then, since \(\gamma(t) \in \mathcal{V}_0\) for all \(t \in [0, 1)\), applying Theorem \ref{teo sobre proy ortog y esfera} to each projection \(\gamma(t)\), we obtain that
		$$
		\gamma(t) = P_{\text{Gr}(A(t))},
		$$
		where \(A(t) = \varphi_0 \begin{pmatrix} \cos|t a^*| \\ (\sinc|t a|) t a^* \end{pmatrix} = (\sinc|t a|) t a^* (\cos|t a^*|)^{-1} = t a^* (\sinc|t a^*|) (\cos|t a^*|)^{-1} = v |t a^*| (\sinc|t a^*|) (\cos|t a^*|)^{-1} = v \tan|t a^*|,\)
		for \(t \in [0, 1)\) and \(v\) is the partial isometry in the polar decomposition of \(a^* = v |a^*|\).
	\end{proof}
			\begin{rem}\label{remark proy sobre graf de no acot en frontera}
				The orthogonal projection onto the graph of any densely defined closed unbounded operator $T$ is in the boundary of the domain of the image of the chart $\varphi_0$ (see \eqref{eq def de fi0 para carta p0}) when we identify the operators $a\in \a$ with their orthogonal projections onto their graphs $P_{\text{Gr}(a)}$ (see Theorem \ref{teo sobre proy ortog y esfera}).
			\end{rem}
	
	\begin{coro}\label{corolario de unicidad de deformacion acotada} 		
		For any unbounded operator $T$ there is a unique bounded deformation $\{T_t\}_{t\in[0,1)}$  (see Definition \ref{def bounded deformation}) such that 
		\begin{enumerate}
			\item $T_0=0$
			\item $t\mapsto \ptilde_t$, $t\in[0,1]$ is a geodesic in $\rcal$
			\item $\ptilde_1=P_{\Gr(T)}$.
		\end{enumerate}
	\end{coro}
	\begin{proof}
		This follows from the properties of the unique minimal geodesic
		$\gamma(t)=P_{\text{Gr}(A(t))}$,  where $A(t)= t a^* (\sinc|t a^*|) (\cos|t a^*|)^{-1}=v \tan|t a^*|\in B(H)$ with $t\in[0,1]$
		from Theorem \ref{teo las geods de PGr(0) a PGr(T) son unicas y son proy sobre graficos}.
	\end{proof}
	\begin{coro}\label{coro deformation length dist}
		The deformation $t\mapsto \{T_t\}_{t\in[0,1]}=\gamma(t)=P_{\text{Gr}(A(t))}$, with $A(t)= t a^* (\sinc|t a^*|) (\cos|t a^*|)^{-1}=v \tan|t a^*|\in B(H)$ is an optimal bounded deformation (see the comments after Definition \ref{def bounded deformation}), that is 
					 $$
					 {\length}|_{t_0}^1 \left(\ptilde_{t_0},\ptilde_1\right)=\dist(\ptilde_{t_0},\ptilde_1)$$
					 	for $t_0\in[0,1)$,  where $\length$ and $\dist$ (distance in the Riemann sphere $\rcal$) are defined for the Finsler metric on $\rcal$ (see Subsection \ref{subsec metrica de Finsler}).
	\end{coro}
%

\subsection{The differential operator}
	
We study here the particular case of an unbounded operator. The conclusions are stated in Theorem \ref{teo props -i d/dx}.

\begin{ejem}[Geodesic between $P_{\gr0}$ and the orthogonal projection onto graph of the differential operator $-i\frac d{dx}$]
	
	Consider the  operator 
	\begin{equation}
		\label{def id/dx}
		-i\frac d{dx}:\mathcal{D}\to L^2[0,1]
	\end{equation}
	given by $f\mapsto -i f'$ for $f:[0,1]\to\mathbb{C}$ with domain 
	\begin{equation}
		\label{def dominio D de i d/dx}
		\mathcal{D}=\{f\in L^2[0,1]: f \text{ is absolutely continuous, }  f'\in L^2[0,1] \text{ and } f(0)=f(1)\}.
	\end{equation}
	This is a known densely defined closed self-adjoint unbounded operator on the Hilbert space $L^2[0,1]$ (see \cite[Example 1.4]{schmudgen}). 
	Denote with 
	$$
	\Gamma=\Gr\left(-i\frac d{dx}\right)=\{(f,-if'): f\in\mathcal{D}\}
	$$ 
	the graph of $-i\frac d{dx}$, which is closed in $L^2[0,1]\times L^2[0,1]$, and $P_\Gamma:L^2[0,1]\times L^2[0,1]\to L^2[0,1]\times L^2[0,1]$ the orthogonal projection onto $\Gamma$.  
	
	Using Lemma \ref{lema Graf T perp} and the fact that $-i \frac d{dx}$ is self-adjoint we have that
	\begin{equation}\label{eq perp del Gr de -id/dx}
		\begin{split}
			\Gamma^\perp&=\left\{\left(-\left(-i \frac d{dx}\right)^*g,g\right): g\in \mathcal{D}\right\}
			=\left\{\left(i   g',g\right): g\in \mathcal{D}\right\}.
		\end{split}
	\end{equation}	
	Theorem \ref{teo las geods de PGr(0) a PGr(T) son unicas y son proy sobre graficos} establishes that there exists a unique geodesic joining $P_{\text{Gr}(0)}$ with $P_\Gamma $. It also can be seen that 
	\begin{equation}\label{H11 y H00}
		\begin{split}
			\mathcal{H}_{11}&:=\ran (P_{\text{Gr}(0)})\cap \ran (P_\Gamma)=[1]\times\{0\}, \\
			\ \text{ and }\ 
			\mathcal{H}_{00}&:=\ker (P_{\text{Gr}(0)})\cap \ker (P_\Gamma)=\{0\}\times[1]
		\end{split}
	\end{equation}
	where $[1]=\{f\in L^2[0,1]: f=\lambda\ 1, \lambda\in \mathbb{C}\}$. And, as in the general case, the matrix block decomposition of $P_{\text{Gr}(0)}$ and $P_\Gamma$ in $\mathcal{H}_{11}\oplus \mathcal{H}_{00}$ is  
	\begin{equation}\label{ec PG y Pgama en H11 + H00}
		P_{\text{Gr}(0)}|_{\mathcal{H}_{11}\oplus \mathcal{H}_{00}}=\left(\begin{smallmatrix}	1&0\\0&0	\end{smallmatrix}\right)=P_\Gamma|_{\mathcal{H}_{11}\oplus \mathcal{H}_{00}}. 
	\end{equation}
	In order to obtain a geodesic between $P_{\text{Gr}(0)}$ and $P_\Gamma$, we need to study $P_\Gamma|_{\mathcal{H}_0}$, the orthogonal projection onto $\Gamma$ when it is restricted to  
	\begin{equation}\label{Hsub0}
		\mathcal{H}_0:= \left(\mathcal{H}_{11}\oplus \mathcal{H}_{00}\right)^\perp	=[1]^\perp \times [1]^\perp. 		
	\end{equation}
	where $[1]^\perp=\{h\in L^2[0,1]:\int_0^1 h(x)\ dx=0\}$. It is clear that $\mathcal{H}_{11}$, $\mathcal{H}_{00}$,  and $\mathcal{H}_{0}$ reduce $P_{\text{Gr}(0)}$ and $P_\Gamma$.
	
	To describe the geodesic that connects $P_{\text{Gr}(0)}$ with $P_\Gamma$, we need to calculate a more specific expression of $P_\Gamma$ restricted to $\mathcal{H}_0=[1]^\perp\times [1]^\perp $.
	Using the  Fourier basis $\{\xi_n\}_{n\in \mathbb{Z}}$ of  $L^2[0,1]$  $\xi_n(x)=e^{i2\pi nx}$, it	follows that 
	$
	\left\{\frac1{\sqrt{1+(2\pi n)^2}}(\xi_n,2\pi n \ \xi_n)\right\}_{n\in \mathbb{Z}\setminus\{0\}}			 
	$ 
	is an orthonormal basis of $[1]^\perp\times [1]^\perp$.
	Then, using blocks in the basis $\{\xi_n\}_{n\in \mathbb{Z}\setminus\{0\}}$ of $[1]^{\perp}$ we obtain that for $(h,k)\in [1]^\perp\times [1]^\perp$
	\begin{equation}\label{ec iddx como diags en base xin}			 				
		\begin{split}
			P_\Gamma|_{[1]^\perp\times [1]^\perp}\left(\begin{smallmatrix}	h\\k\end{smallmatrix}\right)=\left(\begin{smallmatrix}	D_{1}h&D_{2}k\\D_2h&D_3 k\end{smallmatrix}\right)
	\end{split}\end{equation}
	for the following diagonal operators in the  $\{\xi_n\}_{n\in \mathbb{Z}\setminus\{0\}}$ basis 
	\begin{equation}\label{def D1 D2 D3}
		\begin{split}
			D_1	&=\text{diag}\left(\left\{1/\left(1+(2\pi n)^2\right)\right\}_{n\in\mathbb{Z}\setminus\{0\}}\right)
			, \ 
			D_2=
			\text{diag}\left(\left\{{2\pi n}/\left(1+(2\pi n)^2\right)\right\}_{n\in\mathbb{Z}\setminus\{0\}}\right)
			\\
			D_3&=
			\text{diag}\left(\left\{ {(2\pi n)^2}/\left(1+(2\pi n)^2\right)\right\}_{n\in\mathbb{Z}\setminus\{0\}}\right).
		\end{split}
	\end{equation}

	Note that $D_1$ and $D_2$ are positive semidefinite compact operators and $D_3$ is positive definite (invertible) bounded in $[1]^\perp$.
	
	Following ideas from \cite{survey proyectores, halmos2subspaces} and splitting the basis $\{\xi_n\}_{n\in\mathbb Z\setminus \{0\}}$ in $\{\xi_n\}_{n<0}\cup \{\xi_n\}_{n>0}$, we can construct,  the self-adjoint operator $Z_0:[1]^\perp\times [1]^\perp\to[1]^\perp\times [1]^\perp$
	\begin{equation}
		\label{def unitario Z0}
		Z_0=i  \left(
		\begin{smallmatrix}
			0 & 0 & \diag\limits_{n<0}\left\{ -a_n\right\}
			& 0 \\
			0 & 0 & 0 & \diag\limits_{n>0}\left\{ a_n\right\}\\
			\diag\limits_{n<0}\left\{ a_n\right\}  & 0 & 0
			& 0 \\
			0 & \diag\limits_{n>0}\left\{ -a_n\right\}  & 0
			& 0 \\
		\end{smallmatrix}
		\right)
	\end{equation}
	for 
	\begin{equation}
		\label{def a sub n coefs de Z0}
		a_n=\cos^{-1}\left(\frac{1}{\sqrt{4 \pi ^2 n^2+1}}\right)=\left\{
		\begin{array}{rcl} \tan ^{-1}(2\pi n) &\text{, if}& n>0\\ -\tan ^{-1}(2\pi n) &\text{, if}& n<0			
		\end{array}\right..
	\end{equation} 
	
	Observe that  $0<\cos ^{-1}\left(\frac{1}{\sqrt{4 \pi^2 +1}}\right)<a_n=a_{-n}<\pi/2$ for every $n\in \mathbb{N}$ and $\lim_{n\to \pm\infty}a_n=\pi/2$.
	
	Then the unitary $e^{i Z_0}$ satisfies 
	$\label{ec exp iZ0 conjugando al op cero da i d/dx}
	e^{iZ_0}P_{\text{Gr}(0)}|_{[1]^\perp\times [1]^\perp}
	e^{-iZ_0}=\left( \begin{smallmatrix}	D_{1} &D_{2} \\D_2 &D_3  \end{smallmatrix} \right)
	$
	(see \eqref{ec iddx como diags en base xin}, \eqref{def D1 D2 D3}).

	Using the same representation considered in \eqref{def unitario Z0}, since $\lim_{n\to  \pm \infty} \cos ^{-1}\left(\frac{1}{\sqrt{4 \pi^2  n^2+1}}\right)=  \pi/2$ and $0< \cos ^{-1}\left(\frac{1}{\sqrt{4 \pi^2  n^2+1}}\right)<\pi/2$, $\forall n\in\mathbb{Z}$ , we have that $\|Z_0\|=\pi/2$.
	Then we can apply the results from \cite{pr minimality of geod in Grassmann} or \cite[Theorem 5.3]{survey proyectores}: the curve $\delta:[-1,1]\to \mathcal{P}([1]^\perp\oplus [1]^\perp)$
	\begin{equation}\label{ec def de geod delta entre PG y Pgama}
		\delta(t)=e^{it Z_0}\,P_{\text{Gr}(0)}|_{[1]^\perp\oplus [1]^\perp} \,e^{-it Z_0}, \text{ for } t\in[-1,1]
	\end{equation}
	is minimal along its path considering the Finsler metric defined by the operator norm in $[1]^\perp\times [1]^\perp$. 	
	In particular, using our previous computations, this minimal geodesic joins $\delta(0)=P_{\text{Gr}(0)}|_{[1]^\perp\times [1]^\perp}$ with $\delta(1)=P_\Gamma|_{[1]^\perp\oplus [1]^\perp}$.
	Moreover, this implies that (restricted to $[1]^\perp\times [1]^\perp$) the geodesic distance $\text{dist}(P_{\text{Gr}(0)},P_\Gamma)=\pi/2$ and $\|P_{\text{Gr}(0)}-P_\Gamma\|=1$.  
	
	
	Now considering  \eqref{H11 y H00}, \eqref{ec PG y Pgama en H11 + H00}, \eqref{Hsub0}, \eqref{def unitario Z0}, \eqref{ec def de geod delta entre PG y Pgama} and the decomposition $(\mathcal{H}_{11}\oplus \mathcal{H}_{00})\oplus \mathcal{H}_0=L^2([0,1])\times L^2([0,1])$, we can describe the minimal geodesic $\gamma:[-1,1]\to \mathcal{P}( L^2[0,1]\times L^2[0,1])$ 
	\begin{equation}\label{ec geod gama entre PsubG y PsubGama}
		\begin{split}
			\gamma(t)&= \begin{pmatrix}
				\left(\begin{smallmatrix}
					1&0\\0&0
				\end{smallmatrix}\right)&0\\0&\delta(t)
			\end{pmatrix}
			= \begin{pmatrix}
				1&0\\0&e^{it Z_0}
			\end{pmatrix}
			\begin{pmatrix}
				\left(\begin{smallmatrix}
					1&0\\0&0
				\end{smallmatrix}\right)&0\\0&P_{\text{Gr}(0)}|_{[1]^\perp\oplus [1]^\perp}
			\end{pmatrix}					
			\begin{pmatrix}
				1&0\\0&e^{-it Z_0}
			\end{pmatrix}
			\end{split}
	\end{equation}
	with $\gamma(0)=P_{\text{Gr}(0)}$ and $\gamma(1)=P_\Gamma$.
	
	Now observe that the unitary $e^{i t Z_0}$ (see \eqref{def a sub n coefs de Z0})
	in its $2\times 2$ block decomposition, but restricted to  $\{\xi_n\}_{n\in\mathbb{Z}\setminus\{0\}}\times \{\xi_n\}_{n\in\mathbb{Z}\setminus\{0\}}$, can be expressed as $e^{it Z_0}=\left(\begin{smallmatrix}
		A(t)& -B(t)\\ B(t)& A(t)
	\end{smallmatrix}	\right)$, where the diagonal operators $A(t)$ and $B(t)$ are self-adjoint and invertible for $0<t<1$. 
	
	It can also be shown that for $0\leq t<1$ the image of the projection $\delta(t)$ is also the graph of the self-adjoint bounded operator in $[1]^\perp\subset L^2[0,1]$ given by
	\begin{equation}
		B(t)A(t)^{-1}=\diag\limits_{n\in\mathbb{Z}\setminus\{0\}}\{  \tan \left(t \tan ^{-1}(2 \pi  n)\right)\}, \text{ for } 0<t<1.
	\end{equation}
	with $-\tan\left( \frac{t\pi}{2}\right)< \tan\left(t \tan^{-1}(2 \pi x)\right)< \tan\left( \frac{t\pi}{2}\right)$, $\forall x\in \RR$ and norm $\|B(t)A(t)^{-1}\|=\tan\left( \frac{t\pi}{2}\right)$.
	Therefore
	\begin{equation}\label{limite norma -BAala-1}
		\lim_{t\to 1}\|B(t)A(t)^{-1}\|=\lim_{t\to 1}\tan\left( {t\pi}/{2}\right)=+\infty.
	\end{equation}
	Now, considering elements of the whole space $\left(\begin{smallmatrix}	f\\g	\end{smallmatrix}\right)\in 
	L^2[0,1]\times L^2[0,1]= \mathcal H_{11}\oplus\mathcal{H}_{00}\oplus\mathcal{H}_0=
	[1]\times\{0\}\oplus \{0\}\times [1]\oplus[1]^\perp \times [1]^\perp$, we can write
	\begin{equation}\label{eq 1era ecuacion gama}
		\gamma(t)\left(\begin{smallmatrix}f\\g\end{smallmatrix}\right)=
		\left(\begin{smallmatrix}\int_0^1 f\\ 0	\end{smallmatrix} \right)+\left(\begin{smallmatrix}
			A(t)^2& A(t)B(t)\\B(t)A(t)& B(t)^2
		\end{smallmatrix}	\right)
		\left(\begin{smallmatrix}	
			h \\ k 	\end{smallmatrix}\right)=
		\left(\begin{smallmatrix}
			\hat A_d(t)^2 & \ \hat A(t)\hat B(t)\\
			\hat B(t)\hat A(t)&\  \hat B(t)^2
		\end{smallmatrix}	\right)
		\left(\begin{smallmatrix}	
			f \\ g	\end{smallmatrix}\right)
	\end{equation}
	where we denote $\hat A(t)$ and $\hat B(t)$ the corresponding operators extended to $L^2[0,1]$ such that $\hat A(t)(1)=\hat B(t)(1)=0$,  $\hat A_d(t)=\hat A(t)+\diag\{d_j\}_{j\in\mathbb{Z}}$ with $d_j=0$ for $j\ne 0$ and $d_0=1$. We also use that $\hat A_d(t)\hat B(t)=\hat B(t)\hat A_d(t)$ and that $\hat A(t)^2+\diag\{d\}=\hat A_d(t)^2$.
	Also note that $\hat A_d(t):L^2[0,1]\to L^2[0,1]$ is an invertible operator for $0\leq t<1$. 
	
	Hence with the notation used in Section \ref{secc Differential structure of R} and  considering $\xb=\left(\begin{smallmatrix}
		\hat A_d(t) \\ \hat B(t)
	\end{smallmatrix}	\right)$, the $x_1=\hat A_d(t)$ coordinate is invertible for $-1< t<1$, which implies that all the elements $\gamma(t)\in \a^2=B(L^2[0,1])^2$ belong to the chart defined in \eqref{eq def V0 para p0} and \eqref{eq def de fi0 para carta p0}. 
	
	We know from Theorem \ref{teo las geods de PGr(0) a PGr(T) son unicas y son proy sobre graficos} that $\gamma(t)$ are projections onto the graph of an operator for every $t$. In this case it can be proved that, in terms of the Fourier basis, 
	\begin{equation}\label{eq geod op dif como proy sobre graficos}
		\gamma(t)=P_{\Gr(\hat B(t)\hat A_d(t)^{-1})}. 
	\end{equation}
	Therefore, since $\gamma(1)=P_\Gamma$, the entire geodesic $\gamma:[0,1]\to \mathcal{P}( L^2[0,1]\times L^2[0,1])$ is made of self-adjoint orthogonal projections onto graphs of operators. 
	
	Now denote with $\hat D_i$ the diagonal operators such that $\hat D_i(\xi_0)=0$ and $ \hat D_i(\xi_n)=D_i(\xi_n)$, for $n\neq 0$, where $D_i$ are the ones obtained in \eqref{def D1 D2 D3}, for $i=1,2,3$.  Then we have that $\gamma(1)=P_\Gamma=\left(\begin{smallmatrix} \hat D_1+\diag\{d\}&\hat D_2\\ \hat D_2& \hat D_3	\end{smallmatrix}\right)$ expressed in terms of the basis $\left\{\xi_n\right\}_{n\in\mathbb{Z}}\times \left\{\xi_n\right\}_{n\in\mathbb{Z}}$. Since $\gamma(1)_{1,1}=\hat D_1+\diag\{d\}$ is a compact operator then it is not invertible and hence is not in the domain of the chart defined in \eqref{eq def V0 para p0} of Section \ref{secc Differential structure of R}. Moreover, since $\gamma(t)_{1,1}$ is invertible for $0\leq t<1$ and $\gamma(0)=P_{\text{Gr}(0)}=\left(\begin{smallmatrix} 1&0\\ 0&  0	\end{smallmatrix}\right)$, we conclude that $\gamma(1)=P_\Gamma$ lies in the boundary of the principal chart $\mathscr{C}_0$ defined in \eqref{eq def V0 para p0}, a fact that was proven in general in Theorem \ref{teo las geods de PGr(0) a PGr(T) son unicas y son proy sobre graficos} and Remark \ref{remark proy sobre graf de no acot en frontera}. 
	
	We may summarize the above considerations as follows.
	
	\begin{teo}\label{teo props -i d/dx}
		The unique geodesic $\gamma:[0,1]\to \mathcal{P}( L^2[0,1]\oplus L^2[0,1])$ defined in \eqref{eq 1era ecuacion gama} that joins the orthogonal projection $\gamma(0)=P_{\text{Gr}(0)}=P_{L^2[0,1]\oplus\{0\}}$ onto the graph $Gr(0)$ of the zero operator with the orthogonal projection $\gamma(1)=P_\Gamma$ onto the graph $\Gamma$ of the self-adjoint (unbounded, densely defined and closed) differentiation operator $-i\, \frac d{dx}$ (see \eqref{def id/dx}) satisfies the following properties.
		\begin{enumerate}
			\item   For every $t\in (0,1)$, $\gamma(t)$ is the orthogonal projection onto the graph $G_{T(t)}$ of the diagonal self-adjoint bounded operator $T(t)=\hat B(t)\hat A_d(t)^{-1}:L^2[0,1]\to L^2[0,1]$ (see \eqref{eq geod op dif como proy sobre graficos}) that can be written as 
			\begin{equation*}
				T(t)=
				\left(\begin{smallmatrix}
					\diag\limits_{n\in\mathbb{Z}_{<0}}\{  \tan \left(t \tan ^{-1}(2 \pi  n)\right)\} &0&0
					\\ 
					0&0&0\\0&0&\diag\limits_{p\in\mathbb{Z}_{>0}}\{  \tan \left(t \tan ^{-1}(2 \pi  p)\right)\} 
				\end{smallmatrix}\right)
			\end{equation*}
			in terms of blocks determined by the subspaces generated by the respective subsets of the Fourier basis ($\xi_n(x)=e^{i2\pi nx}$, $n\in \mathbb{Z}$) corresponding to $\left\{\xi_j\right\}_{j\in\mathbb{Z}_{<0}}$, $\left\{\xi_0\right\}$  and $\left\{\xi_j\right\}_{j\in\mathbb{Z}_{>0}}$ of $L^2[0,1]$.
			
			\item For $0<t<1$ the operator norm of $T(t)$ is $\|T(t)\|=\tan(t\pi/2)$ and hence 
			$$
			\lim_{t\to 0}\|T(t)\|= 0\text{ and } \lim_{t\to 1}\|T(t)\|=+\infty
			$$ 
			(see \eqref{limite norma -BAala-1} and the properties of $B(t)A(t)^{-1}$).
			\item Every projection $\gamma(t)$ with $0<t<1$, as an element of $\a^2=B(L^2[0,1])^2$,  belongs to the chart defined by \eqref{eq def V0 para p0} and \eqref{eq def de fi0 para carta p0}.  
			\item $\gamma(1)=P_\Gamma$ does not belong to the principal chart $\mathscr{C}_0$ defined in \eqref{eq def V0 para p0} and \eqref{eq def de fi0 para carta p0}; nevertheless, $P_\Gamma$ lies on the boundary of this chart.
			\item $\{T(t)\}_{t\in[0,1)}$ is an optimal bounded deformation of $-i\frac d{dx}$ (see Definition \ref{def bounded deformation}).
		\end{enumerate}
	\end{teo}
	
\end{ejem}

	\subsection{Conjugate parameter values}
	There are infinitely many geodesics  joining two orthogonal projections $\pcero$ and $\ptilde$ in the Grassmann manifold $\Grass(H)$, if an only if $\dim(\ran(\pcero)\cap\ker(\ptilde))=\dim(\ker(\pcero)\cap\ran(\ptilde))\neq 0$ (see \cite{survey proyectores}). 
	Recall that $(T\rcal)_\pcero$ consists of matrices of the form  $X=\begin{psmallmatrix}
		0&a\\a^*&0
	\end{psmallmatrix}$ for $a\in B(H)$. Also, the unique geodesic $\delta$ that satisfies the initial conditions $\delta(0)=\pcero$ and $\dot \delta(0)=X=\begin{psmallmatrix}
		0&a\\a^*&0
	\end{psmallmatrix}$,  is (see \cite[Proposition 2.9]{survey proyectores})
	$$
	\delta(t)=e^{t \tilde{X}} P e^{-t\tilde{X}}, \text{ for } 	\tilde X=\begin{psmallmatrix}		0&-a\\a^*&0		\end{psmallmatrix}.
	$$
	 
	 \begin{teo}\label{teo todas las infinitas geods}
		Let $P_\gr0$ and $Q$ be orthogonal projections such that $\dim(\ran(P_\gr0)\cap\ker(Q))=\dim(\ker(P_\gr0)\cap\ran(Q))\neq 0$. Using the decomposition of $H\times H$ given by $H\times H=\hcal_{11}\oplus\hcal_{00}\oplus  \hcal'\oplus \hcal_0$, with 
		$$
		\hcal'=(\ran(P_\gr0)\cap\ker(Q))\oplus(\ker(P_\gr0)\cap\ran(Q)),
		$$ 
		the geodesics $\gamma:[0,1]\to \rcal$ joining $P_\gr0$ with $Q$ and $\length(\gamma)\leq\pi/2$ are of the form
		\begin{equation}\label{eq todas las geods con |a|<pi/2}
			\gamma_u(t)=\begin{psmallmatrix}
				1 & 0& 0&0\\
				0&0&0&0\\
				0&0&\begin{psmallmatrix} 	\cos^2(t \pi/2)&  \cos(t \pi/2)(\sin(t \pi/2))u\\ (\sin(t \pi/2))(\cos(t \pi/2))u^*& \sin^2(t \pi/2)	\end{psmallmatrix}&0\\
				0 & 0& 0&\delta_0(t)
			\end{psmallmatrix} 
		\end{equation}
 		where: 
 		\begin{itemize}
 			\item  		$u$ is any isometric isomorphism between $\ran(P_\gr0)\cap\ker(Q)$ and $\ker(P_\gr0)\cap\ran(Q)$,
 			\item 		$\dot\gamma|_{\hcal'}(0)=X=\begin{psmallmatrix}
			0&\frac\pi2 u\\\frac\pi2 u^*&0
		\end{psmallmatrix}$, 
		 \item $\delta_0$ is the unique geodesic between the reductions of $P_\gr0$ and $Q$ to $\hcal_0$,
		 \item and $\gamma_u$ has minimal length $\pi/2$.
		\end{itemize}
		\end{teo}
	
	\begin{proof} 
		The multiplicity of these geodesics only appears in  $\hcal'$, which reduces $P_\gr0$ and $Q$ to the expressions $\begin{psmallmatrix}1&0\\0&0\end{psmallmatrix}$ and $\begin{psmallmatrix}0&0\\0&1\end{psmallmatrix}$, respectively (see \cite[Section 3]{survey proyectores}). In what follows we will focus on the geodesics restricted to $\hcal'$.

		Observe that the tangent space at $P_\gr0|_{\hcal'}$ is also formed by co-diagonals $X=\begin{psmallmatrix}		0&a\\a^*&0		\end{psmallmatrix}$ but with $a\in B(\hcal_{01},\hcal_{10})$, and the geodesics starting at $P_\gr0$ are described as 
		$e^{t\tilde{X}}P_\gr0e^{-t\tilde{X}}$ for $\tilde X=\begin{psmallmatrix}		0&-a\\a^*&0		\end{psmallmatrix}$. 
		
		Similarly as we computed in \eqref{eq e a la Xtilde 1} we obtain that	
		\begin{equation}\label{eq unitarios e a la Xtilde}
			e^{\tilde{X}}=\begin{pmatrix}
				\cos|a^*|& -(\sinc|a^*|)a\\(\sinc|a|)a^*&\cos|a|
			\end{pmatrix}
			=\cos|\tilde X|+\left(\sinc|\tilde X|\right)\tilde X,
		\end{equation}
		because $|\tilde X|=\begin{psmallmatrix}
			|a^*|&0\\0&|a|
		\end{psmallmatrix}$.
		Then using that $e^{\tilde{X}}$ is a unitary operator and must satisfy $e^{\tilde{X}}\begin{psmallmatrix}
			1& 0\\0&0\end{psmallmatrix}e^{-\tilde{X}}=\begin{psmallmatrix}
			0& 0\\0&1\end{psmallmatrix}$ it can be proved that $\cos|a^*|=\cos|a|=0$ and hence $(\sin|a^*|)^2=(\sin|a|)^2=1$.
			Thus we obtain that $|a|=\sum_{j=1}^m (n_j \pi+\pi/2) p_j$ with $n_j\in \mathbb{N}\cup \{0\}$
		for $p_j$ spectral projections of $|a|$ that satisfy $\sum_{j=1}^n  p_j=1$. Considering that $\Xtilde$ must satisfy $\|\Xtilde\|=\|a\|=\|a^*\|\leq \pi/2$ ($\length(\gamma)\leq \pi/2$), we have that	$|a|=\frac{\pi}2$ and then $a=\frac{\pi}2u$, with $u:\hcal_{01}\to\hcal_{10}$ a unitary isomorphism.
		And then $a^* =\frac\pi2 u^*$.
		
		Therefore all the possible $\tilde X$ are of the form 
		$\tilde{X}=
		\begin{pmatrix}
			0&-\frac\pi2u\\ \frac\pi2 u^*&0
		\end{pmatrix}
		$.
		And all the unitaries $e^{t\tilde X}$ are (see \eqref{eq unitarios e a la Xtilde})
		$$
		e^{t \tilde{X}}=\begin{pmatrix}
			\cos|t \frac\pi2|& -(\sinc|t\frac\pi2|)t \frac\pi2 u\\(\sinc|t \frac\pi2|)t \frac\pi2 u^*&\cos|t \frac\pi2|
		\end{pmatrix}
		=\begin{pmatrix}
			\cos(t \frac\pi2)& -\sin(t\frac\pi2)u\\ \sin(t \frac\pi2)u^*&\cos(t\frac\pi2)
		\end{pmatrix}.
		$$
					
			Then considering the decomposition of $H\oplus H=\hcal_{00}\oplus \hcal_{11}\oplus \hcal'\oplus \hcal_0$ with $\hcal'=\hcal_{10}\oplus\hcal_{01}$ (see \cite[Section 3]{survey proyectores}) all the geodesics between the projections $P_\gr0$ and $Q$ can be parameterized more explicitly as
			\begin{equation*}\label{eq todas las infinitas geos entre P y Q}
				\begin{split}
					\gamma(t)&=
					\begin{pmatrix} \begin{psmallmatrix}	1 & 0\\ 0&0 \end{psmallmatrix} &0&0\\
						0&e^{t\tilde{X}}\begin{psmallmatrix}1&0\\0&0	\end{psmallmatrix} e^{-t\tilde{X}}&0\\
						0&0&\delta_0(t)
					\end{pmatrix}
					=\begin{psmallmatrix}
						1 & 0& 0&0\\
						0&0&0&0\\
						0&0& \begin{psmallmatrix} 
							\cos^2(t \frac\pi2)& \cos(t \frac\pi2)\sin(t \frac\pi2)u\\
							\sin(t \frac\pi2)u^*\cos(t \frac\pi2)&\sin^2(t \frac\pi2)
								\end{psmallmatrix}&0\\
						0 & 0& 0&\delta_0(t)
					\end{psmallmatrix} 
				\end{split}
			\end{equation*}
			which is the expression for the geodesics $\gamma_u$ stated in \eqref{eq todas las geods con |a|<pi/2}.
			
			The minimality condition of the geodesics $\gamma$ when $t\in [0,1]$ and $\|a\|\leq \frac\pi2$ follows from \cite[Theorem 5.3 and Corollary 5.5]{survey proyectores}.
		\end{proof}

		Recall that a classical Jacobi field is a field on a fixed geodesic $\gamma$ that can be obtained differentiating a family of perturbations of $\gamma$ by geodesics that start and end at the same points as $\gamma$.
		\begin{defi}
			Given the geodesic $\gamma(t)$, $t\in[0,1]$ a parameter value $t_0\in [0,1]$ is called \textbf{conjugate of $0$ along $\gamma$} if there exists a non trivial Jacobi field that vanishes at $0$ and at $t_0$. In this case the \textbf{index of $t_0$} is the dimension of the space of Jacobi fields that vanish at $0$ and at $t_0$. The parameter $t_0$ is called \textbf{conjugate} if this index is greater than zero. 
		\end{defi} 
		The following is an example of conjugate values in $\rcal$ (see \ref{rem PinvGr(T) esta en R}) involving Fredholm operators.	
		\begin{teo}
			Let $F$ be a Fredholm operator of index zero and call $n=\dim(\ker(T))=\dim(\ran(T)^\perp)> 0$. Then, $1$ is a conjugate parameter of $0$ for the geodesic defined in 
			\eqref{eq todas las geods con |a|<pi/2} for $u=1$, with $t\in[0,1]$, connecting the orthogonal projection onto the graph of the null operator $P_{\text{Gr}(0)}$ and the orthogonal projection onto the inverse graph of $P_{{\invGr(T)}}$.
			
			Moreover, the index of this conjugate parameter has dimension $n^2$.
%
		\end{teo}
		\begin{proof}
			From Lemma \ref{lema Graf T perp} we can state that that $Gr(T)^\perp=\{(-T^*x,x):x\in \text{Dom}(T^*)\}=\text{invGr}(-T^*)$. And since $\ran (P_{\text{Gr}(0)})=H\oplus \{0\}$ and $\ker (P_{\text{Gr}(0)})=\{0\}\times H$, we obtain that
			\begin{equation}\label{eq H10 y H01 para cograf} \begin{split} \mathcal{H}_{10} &= \ran (P_{\text{Gr}(0)}) \cap \ker (P_{\invGr(T)}) = \ker (T^*) \oplus \{0\} = \ran (T)^\perp \oplus \{0\}, \\ \mathcal{H}_{01} &= \ker (P_{\text{Gr}(0)}) \cap \ran (P_{\invGr(T)}) = \{0\} \oplus \ker (T). \end{split} \end{equation}
			
			Therefore the condition $\dim(\ker(T))=\dim(\ran(T)^\perp)=n> 0$ 
			implies that there exist infinite geodesics joining $P_{\text{Gr}(0)}$ with $P_{{\invGr(T)}}$ (see \cite{survey proyectores}).

			Then, we can use Theorem \ref{teo todas las infinitas geods} and the expression of the geodesics $\gamma_u$ from \eqref{eq todas las geods con |a|<pi/2}. We will differentiate curves of geodesics using the parameter $s$ that describes unitaries $u(s)$ with fixed $t$. Observe that the only part that changes is in the $\hcal'=\hcal_{10}\oplus\hcal_{01}$ space given by 
			$$
			{{\gamma_{u}}|}_{\hcal'}(t) =	\begin{psmallmatrix} 	\cos^2(t \pi/2)& -\cos(t \pi/2)(\sin(t\pi/2))u\\ (\sin(t \pi/2))(\cos(t \pi/2))u^*&-\sin^2(t \pi/2)	\end{psmallmatrix}\ \text{ for } t\geq 0
			$$
			We will construct a Jacobi field along the fixed geodesic $\gamma_1$ which is the case when $u=1$. In this case we can consider the Jacobi field obtained after differentiating the geodesics perturbed by unitary curves $u(s)$ that depend on the parameter $s$ close to $s=0$ with $u(0)=1$. Then, differentiating ${{\gamma_{u}}|}_{\hcal'}(t)$ respect to $s$ we have
			$$
			\mathcal{J}_{\dot{u}}(t)=\frac{\partial}{\partial s}\Big|_{s=0}\left( {{\gamma_{u(s)}}|}_{\hcal'}(t) \right)
			=\begin{pmatrix} 	0 & -\cos(t \pi/2)(\sin(t\pi/2))\dot u(s)\\ (\sin(t \pi/2))(\cos(t \pi/2)) \dot u^*(s)&0	\end{pmatrix} 
			$$
			for $t\in[0,1]$. Note that the derivatives $\dot u(s)$ belong to the space of anti self-adjoint operators (elements of the Lie algebra of the unitary group) that has real dimension $n^2$ since $u:\hcal_{01}\to\hcal_{10}$ (each of dimension $n$).
	\end{proof}
	
\subsection{Density (and non density) of the geodesic neighborhoods}

Let us briefly examine examples of algebras where $\{\qtilde\in P_2(\a): \|\qtilde-\pcero\|<1\}$ is dense in the orbit of $p_0$, and examples where it is not. The first example includes the case of finite matrices.
\begin{ejem}\label{finito}
	Let $\a$ be a finite von Neumann factor, with (unique) normal, faithful and normalized trace $\tau$. Then $M_2(\a)$ is also a a finite von Neumann factor with trace ${\bf Tr}\begin{psmallmatrix} a & b \\ c & d\end{psmallmatrix} =\frac12\tau(a+d)$. In \cite{geodesics von neumann} it was shown that any pair of projections in a finite factor, in the same connected component (i.e., in the same unitary orbit, or equivalently, with equal trace) can be joined with a minimal geodesic. Pick as usual $\pcero\in M_2(\a)$, $\ptilde_0= \begin{psmallmatrix} 1 & 0 \\ 0 & 0 \end{psmallmatrix}$. Let $\qtilde$ be a projection in the orbit of $\pcero$, and $\gamma(t)=e^{itX}\pcero e^{-itX}$ a geodesic with $\gamma(1)=\qtilde$, with $X^*=X$ $\pcero$-co-diagonal and $\|X\|\le\pi/2$. It is known that \cite{pr minimality of geod in Grassmann}
	$$
	\|\gamma(t)-\gamma(s)\|=\sin\left(|t-s|\|x\|\right).
	$$
	Therefore, given $\epsilon>0$, we can choose $t_0<1$ such that $\qtilde_0=:\gamma(t_0)$ satisfies
	$$
	\|\qtilde_0-\qtilde\|=\|\gamma(t_0)-\gamma(1)\|=\sin\left((1-t_0)\|X\|\right)<\epsilon.
	$$ 
	Clearly also $\|\pcero-\qtilde_0\|=\|\gamma(0)-\gamma(t_0)\|<1$. That is, $\{\qtilde\in P_2(\a): \|\qtilde-\pcero\|<1\}$ is dense in $\mathcal{R}$, the orbit of $\pcero$.  
\end{ejem}
The next example shows that this is no longer the case if $\a=\mathcal{B}(H)$, for $H$ infinite dimensional. To present the specific subspaces, first we need to recall results on the theory of common complements of pairs of subspaces, as presented by M. Lauzon and S. Treil in \cite{treil}, and continued by J. Giol \cite{giol}.
\begin{rem}\label{complemento comun}
	In \cite{treil}, necessary and sufficient conditions were given, in order that a pair of closed subspaces $\mathcal{S}, \mathcal{T}$ of an infinite dimensional Hilbert space $\mathcal{L}$ do (or do not) have a common complement, i.e., that there exists (or not) a closed subspace $\mathcal{Z}\subset \mathcal{L}$ such that $\mathcal{S}\dot{+}\mathcal{Z}=\mathcal{L}$ and $\mathcal{T}\dot{+}\mathcal{Z}=\mathcal{L}$, where the symbol $\dot{+}$ stand for direct non necessarily orthogonal sum. For instance, in \cite{treil} it was shown that $\mathcal{S}, \mathcal{T}\subset \mathcal{L}$ do not have a common complement if and only if $\dim \mathcal{S}\cap\mathcal{T}^\perp \ne  \dim \mathcal{S}^\perp\cap\mathcal{T}$ and 
$$
1_{\mathcal{S}}-G^*G:\mathcal{S}\to\mathcal{S} \hbox{ is compact when restricted to } N(G)^\perp,
$$
where $G:=P_{\mathcal{T}}\big|_\mathcal{S}:\mathcal{S}\to\mathcal{T}$. Here $N(G)=\mathcal{S}\cap\mathcal{T}^\perp$.
	
	Later on J. Giol \cite{giol} proved that $\mathcal{S}$ and $\mathcal{T}$ do have a common complement if and only if there exists an intermediate orthogonal projection $Q$ such that $\|P_\mathcal{S}-Q\|<1$ and $\|Q-P_\mathcal{T}\|<1$.
\end{rem}
Building on these facts, it is easy to see that a pair of subspaces $\mathcal{S}$, $\mathcal{T}$ of $\mathcal{L}$, with infinite and co-infinite dimension, and  without a common complement, provide an example where $\{P\in P(\mathcal{L}): \|P_\mathcal{S}-P\|<1\}$ is not dense in the unitary orbit of $P_\mathcal{S}$: indeed, note that in this case  $$
\{P\in P(\mathcal{L}): \|P_\mathcal{S}-P\|<1\}\cap\{Q\in P(\mathcal{L}): \|P_\mathcal{T}-Q\|<1\}=\emptyset.
$$
Clearly, an element $Q$ in this intersection would provide an intermediate projection with $\|P_{\mathcal S}-Q\|<1$ and $\|P_{\mathcal T}-Q\|<1$, and this would imply, by Giol's result, that ${\mathcal S}$ and ${\mathcal T}$ have a common complement. 

Also it is clear how to adapt this example to our situation (where one of the subspaces is $H\times\{0\}$). Pick a unitary isomorphism $U:{\mathcal L}\to H\times H$ which maps ${\mathcal S}$ onto $H\times\{0\}$. This is done by choosing orthonormal bases of ${\mathcal S}$ and $H\times\{0\}$, and completing them to orthonormal bases of ${\cal L}$ and $H\times H$, respectively, and is possible because ${\mathcal S}$  has infinite and co-infinite dimension. Since ${\mathcal S}$ and ${\mathcal T}$ do not have a common complement in ${\mathcal L}$, it is clear that $H\times\{0\}=U{\mathcal S}$ and $U{\mathcal T}$ do not have  common complement in $H\times H$. 

Therefore $\{P\in\mathcal{R}: \|\pcero-P\|<1\}$ is not dense in $\mathcal{R}$ in this case.



\begin{ejem}\label{infinito}

This example was discussed in \cite{incertidumbre} in connection with existence and non existence of geodesics between subspaces, and it is related to the so called Uncertainty Principle in Harmonic Analysis.

	\bigskip

	Let $I,J\subset\mathbb{R}^n$ be Lebesgue measurable subsets with finite positive measure. Consider 
	$$
	{\bf S}_I=\{f\in L^2(\mathbb{R}^n): \text{supp}(f)\subset I\}, \ {\bf T}_J=\{g\in L^2(\mathbb{R}^n): \text{supp}(\hat{g})\subset J\},
	$$
	where $\text{supp}$ stands for the (essential) support, and $\hat{g}$ is the Fourier-Plancherel  transform of $g$. Put $\mathcal{S}={\bf S}_I$ and $\mathcal{T}={\bf T}_J^\perp$. We claim that 
 $\mathcal{S}$ and $\mathcal{T}$ do not have a common complement.

Indeed, it is known that (see \cite{lenard} or the survey article \cite{follandsitaram})
$$
\mathcal{S}\cap\mathcal{T}^\perp={\bf S}_I\cap{\bf T}_J=\{0\} \ \hbox{ and } \ \mathcal{S}^\perp\cap\mathcal{T}={\bf S}_I^\perp\cap{\bf T}_J^\perp \ \hbox{ is infinite dimensional}.
$$
Also, it is known  that $P_{{\bf S}_I} P_{{\bf T}_J}$ is compact (see \cite{follandsitaram}). This clearly means that 
$$
P_\mathcal{S} -P_\mathcal{S}P_\mathcal{T}P_\mathcal{S}=P_\mathcal{S}P_\mathcal{T}^\perp P_\mathcal{S}=P_{{\bf S}_I} P_{{\bf T}_J}P_{{\bf S}_I}
$$
is compact, i.e., $1_\mathcal{S}-G^*G$ is compact  in  the whole $\mathcal{S}$ (here $N(G)=\mathcal{S}\cap\mathcal{T}^\perp=\{0\}$). Therefore, by the result of Lauzon and Treil \cite{treil} transcribed in Remark \ref{complemento comun}, $\mathcal{S}$ and $\mathcal{T}$ do not have a common. complement
\end{ejem}
\begin{rem}
	The Example \ref{infinito}  tells us that the classical Hopf-Rinow Theorem is not valid when $\mathcal{A}=\mathcal{B}(H)$ for $H$ infinite dimensional. There are points in $\mathcal{R}$ which cannot be reached by a geodesic starting at $\pcero$, not even approximated by points in the range of the exponential based at $\pcero$. Moreover, elaborating on this example, it  also shows that there exist in $\mathcal{R}$ infinitely many disjoint open subsets,  which are ranges of the exponential map at different points in $\mathcal{R}$.
	
	Example \ref{finito} suggests that density of the range of the exponential at $\pcero$ requires some sort of finiteness (for instance, that the algebra is finite, as shown in this example).
	
\end{rem}

\end{document}